\documentclass[12pt,a4paper]{amsart}
\usepackage[utf8]{inputenc}
\usepackage[in,headings]{fullpage}
\usepackage{amssymb,amscd,amsbsy,amsmath, amsthm, amsfonts, amsrefs}
\usepackage[OT2,OT1]{fontenc}
\usepackage{color}   
\usepackage[colorlinks=true,linkcolor=blue,citecolor=blue,urlcolor=blue]{hyperref}
\usepackage{fancyvrb}

\newtheorem{theorem}{Theorem}

\newtheorem{proposition}[theorem]{Proposition}
\newtheorem{lemma}[theorem]{Lemma}

\theoremstyle{definition}

\theoremstyle{remark}

\newfont{\cmbsy}{cmbsy10}
\newfont{\cmmib}{cmmib10}
\newcommand{\Orden}{\mathop{\hbox{\cmbsy O}}\nolimits}

\def\Turing{\Orden^*}
\def\Z{\mathbf{Z}}
\def\N{\mathbf{N}}
\def\Q{\mathbf{Q}}
\def\R{\mathbf{R}}
\def\C{\mathbf{C}}

\renewcommand{\Re}{\operatorname{Re}}
\renewcommand{\Im}{\operatorname{Im}}
\newcommand{\Rzeta}{\mathop{\mathcal R }\nolimits}

\def\pmm{\mathrm{\pm\,}}
\def\bc{\mathrm{bc}}
\def\({($\,$}
\def\){$\,$)}
\def\sqrtt{\mathrm{sqrt}}
\def\ba{\vbox{\hrule width 6pt height 0.6pt depth 0.2pt }}  

\begin{document}

\title[High Precision Computation of zeta]
{High Precision Computation of Riemann's Zeta Function by the Riemann-Siegel Formula, II.}
\author[Arias de Reyna]{J. Arias de Reyna}
\address{Facultad de Matemáticas \&\ IMUS\\
Universidad de Sevilla\\
c. Tarfia, sn\\
41012-Sevilla \\
Spain} 
\email{arias@us.es}

\date{\today, \texttt{88-RiemannSiegel-v3.tex}}

\begin{abstract}
(This is only a first preliminary version, any suggestions about it will be welcome.)
In this paper it is shown how to compute Riemann's zeta function $\zeta(s)$ (and Riemann-Siegel $Z(t)$) at any point $s\in\C$ with a prescribed error  $\varepsilon$ applying the Riemann-Siegel formula as described in my paper  \emph{High Precision \dots I},  Math of Comp. \textbf{80} (2011) 995--1009.

This includes the study of how many terms to compute and to what precision
to get the desired result. All possible errors are considered, even those 
inherent to the use of floating point representation of the numbers.

The result has been used to implement the computation. The programs have been included in \texttt{mpmath}, a public library in Python for the computation of special functions. Hence they are included also in Sage.
\end{abstract}
\maketitle
\tableofcontents

\setcounter{section}{0}
\section{Introduction.}

In \cite{A86} we give explicit bounds for the rest of the Riemann Siegel expansion. The Riemann Siegel expansion allow us to compute $\zeta(s)$ to a given precision (almost any point and almost any precision). I set out to implement this calculation. Given $s\ne 1$ and  $\epsilon>0$ the task is to compute $a$ such that $|\zeta(s)-a|<\varepsilon$. This is a very complicated project.  We have to answer many questions:  Is this calculation possible? How many terms of the Riemann Siegel development will be necessary? How can we calculate these terms? How to take into account all the errors in the necessary operations?

Knuth \cite{Knuth}*{Section 4.2.2 Accuracy of Floating Point Arithmetic, p.~229} observes: \emph{Many serious mathematicians have attempted to analyze a sequence of floating point operations rigorously, but have found the task so formidable that they have tried to be content with plausibility arguments instead.}
But in modern times,  libraries for arbitrary-precision floating-point arithmetic  have appeared such as \texttt{mpmath} \cite{FJ}. At any moment in the computation we can increase the number of digits in which we operate, so that the result has the desired approximation. 

This is not the first implementation. Notable is the one by Jerry B. Keiper in Mathematica. This computes zeta at any point and to a given precision, but it is not documented. Comparison with our implementation shows that it is very reliable. In fact we use it as a check of our computations. Many implementations are based on the explicit bounds given by W. Gabcke, but they are limite to the critical line.  

This paper was written for my use in the implementation for \texttt{mpmath}. This is free and open source and included in Sage \cite{sage}. Its publication would require many checks, and proofs that are not included. For example the Technical lemmas in Section 3.10 are only checked by means of convincing plots. If not a human readable proof I expect they will be able to be proved by means of the Maximal slope principle as described in \cite{A94}.

I would like to add that for the first time I proof the limits of the Riemann Siegel development. For example it is shown that for $t\in\R$ we can compute $\zeta(\frac12+it)$ with any error greater than $e^{-2t/\pi}$, see Section \ref{S:3.3}.

At that time (October 2009) it was clear that the implementation gave correct values. But at the same time I knew that this work was not ready for publication.  Its subject, its length, that I had not sufficiently tested each formula,  and that some inequalities were not proven (although as I said I do not think these proofs are a real problem), all this made me abandon the purpose of publishing it, its objective, was completed with the publication of my implementation for the computation of zeta. 

Nevertheless, I put references in my code to this paper. Since this code is open some people ask me for a copy of this paper. So, without any revision I put it on arXiv now. I would be grateful to anyone who has any comments to communicate them to me in order to improve this paper. 

I do not want to end this introduction without mentioning the problem detected in the usual implementation of the product of complex numbers. See Section \ref{S:2.6} where I indicate what a good definition should fulfill.

\section{Multiple Precision Floating Point Arithmetic.}

We will  describe  an idealized system of Multiple Precision Floating
Point  Arithmetic \(MPFP\). In practice the program  described here
has been implemented in Python using the \texttt{mpmath} library for MPFP.

\subsection{Representable numbers.}
\emph{Dyadic numbers} are those of the form $\pmm m\cdot 2^e$ where
$m\in\N$ and $e\in\Z$.  \emph{Representable numbers} are dyadic
numbers and $0$.

We assume that the system has an internal representation for each 
\emph{representable}
number. This is clearly not true, but we will assume that in practice
the numbers appearing in a given computation are implementable
on the computer. Therefore we do not consider overflow on exponents
and/or mantissa.

Any dyadic number $x=\pmm A\cdot 2^E$ can be put in a unique way in
a \emph{normalized} form $\pmm m\cdot 2^e$ where $m$ is odd. If
$2^{d-1}\le m<2^{d}$, the dyadic expansion of $m$ has $d$ digits and
we say that $d=:\bc(x)$ (bit count of $x$). So
\begin{equation}
m< 2^{\,\bc(x)}\le 2m.
\end{equation}

\subsection{Rounding.}
Given a dyadic number $x$  and a natural number $d$ we define the
\emph{rounded version} of $x$, denoted by $\circ(x,d)$ \(or $\circ_d(x)$\):

To define it we observe that given $m$ and $2^f$ with $f\in\N$,  by a
modified division algorithm, there exists a representation
\begin{equation}\label{semiquotient}
m=q\cdot 2^{f}\pmm r,\qquad 0\le r\le 2^{f-1}.
\end{equation}
This representation is unique, except when $r=2^{f-1}$ in
which case there are two. To get a unique representation we ``round
to even'', that is we take the only representation with $q$ even.

Then we can define the rounded version $x'=\circ(x,d)$. We define $x'=x$ if
 $\bc(x)\le d$. In the other case we assume that
$x=\pmm m \cdot 2^e$ is the normalized representation of $x$,  and
take the unique representation \eqref{semiquotient} with
\begin{equation}
m=q\cdot 2^{\bc(x)-d}\;\pmm r\;\quad\text{and}\quad 0\le r\le 2^{\bc(x)-d-1}.
\end{equation}
In this case  we define
\begin{equation}
\circ(x,d)=\pmm q \cdot 2^{\bc(x)-d+e}
\end{equation}
where the sign is the same as that of $x$.

When $x=0$ we define $\circ(x,d)=0$.

We always have the important relation
\begin{equation}\label{boundround}
\circ(x,d)=x(1+\eta_1)\qquad \text{with}\quad |\eta_1|\le 2^{-d}.
\end{equation}
When either $x=0$ or $\bc(x)\le d$  this is trivial. In the other case we
have
\begin{multline*}
|x-\circ(x.d)|=| m\cdot 2^e-q \cdot 2^{\bc(x)-d+e}|= \\|(q\cdot
2^{\bc(x)-d}\pm r)\cdot 2^e-q \cdot 2^{\bc(x)-d+e}|=|r\cdot
2^e|\le\\
\le 2^{\bc(x)-d-1} 2^e\le 2 m\cdot 2^{e-d-1}=|x| 2^{-d}
\end{multline*}
which is equivalent to \eqref{boundround}.
\bigskip

The ``round to even'' of $x$ to $d$ bits  $\circ(x,d)$ can be  seen as the
truncation of the dyadic expansion of $x$. If $x=\pmm m \cdot 2^e$
and we expand $m$ in base $2$
\begin{displaymath}
m=\varepsilon_1\varepsilon_2\cdots
\varepsilon_d\varepsilon_{d+1}\cdots
\varepsilon_{\bc(x)}
\end{displaymath}
then
\begin{displaymath}
q=\varepsilon_1\varepsilon_2\cdots \varepsilon_d \quad
\text{or}\quad q=\varepsilon_1\varepsilon_2\cdots \varepsilon_d+1.
\end{displaymath}
The $1$ has to be included if $\varepsilon_{d+1}=1$, except in the
case that all the other digits $\varepsilon_j$ for $j\ge d+1$ are
equal to $0$ and $\varepsilon_d=0$.
\medskip

It is easy to extend the definition of $\circ(x,d)$ to the case of
an arbitrary real  $x$, not necessarily a dyadic number. The
relation \eqref{boundround} is also true in this more general setting.

\subsection{Working Precision and Elementary Operations.}
In an MPFP system  there is a variable \(\texttt{mp.prec} in
mpmath\) that we can set equal to any natural number $\ge1$. The value
$d$ of this variable is called \emph{working precision}.  The system
works as a floating point system with precision equal to $d$.

When we operate with representable numbers  $x$ and $y$  the result
is as if we operate with $\circ_d(x)$ and $\circ_d(y)$. So if we put
\begin{displaymath}
\texttt{z = x + y},\quad\texttt{z = x - y},\quad \texttt{z = x *
y},\quad\texttt{z = x / y}
\end{displaymath}
the computer will put in  $z$ the representable numbers
\begin{multline*}
z=\circ_d(\circ_d(x)+\circ_d(y)),\quad
z=\circ_d(\circ_d(x)-\circ_d(y)),\\
z=\circ_d(\circ_d(x)\cdot\circ_d(y)),\quad
z=\circ_d(\circ_d(x)/\circ_d(y)).
\end{multline*}
This is in accord with the IEEE standard requirement that the result
of addition, subtraction, multiplication and division must be
computed exactly, and then rounded to the nearest floating-point number
\(using round to even\).

\(Observe that in the case of division, in general,
$\circ_d(x)/\circ_d(y)$ is not dyadic, but we have defined
$\circ_d(x)$ even for any real number $x$.\)

The elementary operations are performed in two steps: first the numbers
$x$ and $y$ are rounded and  then the operation is performed. Since
the first rounding step introduces truncation errors we shall
introduce in our algorithms \emph{fictitious steps} to take account
of these truncation errors.

So we shall put instead of \texttt{z = x + y}
\begin{equation}\label{fictitious}
\begin{split}
&\texttt{mp.prec=d}\\
&\texttt{x1 = round(x,d); y1 = round(y,d)}\\
&\texttt{z = x1 + y1 }.
\end{split}
\end{equation}
We will also use \texttt{x1=round(x)} understanding that, if not
explicit, $d$ is the actual value of the working precision.

We shall say that $x$ is \emph{rounded} if $\bc(x)\le d$ where $d$ is the working precision.
So $x$ is rounded if and only if $x=\circ_d(x)$.  

By \eqref{boundround} IEEE standard requirement implies  the following Proposition.
\begin{proposition}\label{IEEEstandard}
Let $x$ and $y$ be rounded real numbers. Then  there exist real
numbers $\eta_1$ such that
\begin{align*}
&x\oplus y=(x+y)(1+\eta_1), \quad & |\eta_1|\le 2^{-d}\\
&x\ominus y=(x-y)(1+\eta_1),\quad & \quad |\eta_1|\le 2^{-d}\\
&x\otimes y=xy(1+\eta_1), \quad & |\eta_1|\le 2^{-d},
\end{align*}
and if $y\ne0$
\begin{align}
&x\oslash y=\frac{x}{y}(1+\eta_1), \quad & |\eta_1|\le 2^{-d},
\end{align}
\end{proposition}

\subsection{Notation  $\eta_k$.}
Due to \eqref{boundround} in the study of the computation errors there
will  frequently appear products of the form
\begin{equation}
1+\eta_k=\prod_{j=1}^k (1+\xi_j)\qquad \text{where for each $j$,
}\quad |\xi_j|\le 2^{-d}.
\end{equation}
We will  always use the notation $1+\eta_k$ to denote such a
product of $k$ factors. So, two instances of $\eta_k$ do not represent the same
number, even when they are part of the same formula. This poses no
problem since we are only interested in the bound of the absolute
value of these $\eta$ numbers.

So, with this notation we have the following very useful relation
\begin{equation}
(1+\eta_r)(1+\eta_s)=(1+\eta_{r+s}).
\end{equation}

The following Lemma ($\,$an expanded version of one given by Wilkinson
\cite{Wilkinson}*{p.~19}$\,$)  gives us  an adequate bound for
$\eta_r$.

\begin{lemma}\label{boundExp}
Let $r$ and $d$ be  natural numbers such that $r\cdot2^{-d}\le0.1$.
Assume that for some complex numbers $\xi_j$ we have
\begin{displaymath}
1+\eta_r=\prod_{j=1}^r (1+\xi_j)\qquad \text{ with\quad $|\xi_j|\le
2^{-d}$ for } 1\le j\le r.
\end{displaymath}
Then we will have
\begin{equation}
|\eta_r|\le 1.06 \cdot r\cdot 2^{-d}.
\end{equation}
\end{lemma}

\begin{proof}
First we assume that $\xi_j>0$ for each $j$. Then  $\eta_r>0$.
Putting $x=2^{-d}$, since $\binom{r}{j}\le\frac{r^j}{j!}$,  we have
\begin{multline*}
\eta_r\le (1+x)^r-1=rx+\sum_{j=2}^r\binom{r}{j} x^j\le
rx+rx\Bigl(\frac{rx}{2!}+ \frac{r^2x^2}{3!}+\cdots
\frac{r^{n-1}x^{n-1}}{n!}\Bigr) \le\\\le
rx\Bigl(1+\frac{e^{rx}-1}{2}\Bigr)\le
rx\Bigl(1+\frac{e^{0.1}-1}{2}\Bigr)\le rx \cdot 1.0525855 \le 1.06\,
r\, x.
\end{multline*}
Now in the general case we have
\begin{displaymath}
\eta_r=\prod_{j=1}^r (1+\xi_j)-1=\sum_{J\subset\{1,2,\dots,r\}}
\xi_J
\end{displaymath}
where the sum extends over all non-empty subsets
$J\subset\{1,2,\dots,r\}$, and we define $\xi_J=\prod_{j\in J}
\xi_j$.

Therefore
\begin{displaymath}
|\eta_r|\le\sum_{J\subset\{1,2,\dots,r\}}
|\xi_J|=\prod_{j=1}^r(1+|\xi_j|)-1.
\end{displaymath}
Since $|\xi_j|\le 2^{-d}$ for all $j$  the conditions of
 the previous case are satisfied so that as before we get
\begin{displaymath}
|\eta_r|\le 1.06\,r\,x.\qedhere
\end{displaymath}
\end{proof}

For example, with   the notations here introduced, the code in  \eqref{fictitious}
corresponds with the following bounds
\begin{displaymath}
x_1 =x(1+\eta_1),\quad y_1=y(1+\eta_1),\quad z=(x_1+y_1)(1+\eta_1)
\end{displaymath}
so that $z=x(1+\eta_2)+y(1+\eta_2)$.  Observe that this is not
equal to $z=(x+y)(1+\eta_2)$, since $\eta_2$ is a different constant
in each case.

\begin{lemma}\label{lemma2}
A factor of the form $(1+\eta_1)^{-1}$ can be written as $(1+\eta_2)$.
\end{lemma}

\begin{proof}
Define $\xi$ by $1+\xi=(1+\eta_1)^{-1/2}$ \(the value given by Newton
series\). We must show that $|\xi|\le 2^{-d}$.

By definition
\begin{displaymath}
\xi=\sum_{n=1}^\infty \binom{-1/2}{n} \eta_1^n
\end{displaymath}
so that
\begin{displaymath}
|\xi|\le\sum_{n=1}^\infty\Bigl|\binom{-1/2}{n}\Bigr|2^{-dn}=
\sum_{n=1}^\infty(-1)^n \binom{-1/2}{n}2^{-dn}=(1-2^{-d})^{-1/2}-1
<2^{-d}.\qedhere
\end{displaymath}

\end{proof}

\begin{lemma}\label{lemma3}
Let $d'>d\ge1$. Let $A=(1+\eta'_n)$ be a factor associated to the working precision
$d'$. If we change the working  precision to $d$, then $A^{-1}=(1+\eta_n)$ where $\eta_n$ is relative to the new working precision.
\end{lemma}

\begin{proof}
We know that $A=(1+\eta'_1)\cdots (1+\eta'_1)$. It is clear that we only need to consider the
case $n=1$.  Let us define $\xi$ such that
$(1+\xi)=(1+\eta'_1)^{-1}$. Then
\begin{displaymath}
|\xi|=|(1+\eta'_1)^{-1}-1|=\Bigl|\frac{\eta'_1}{1+\eta'_1}\Bigr|\le \frac{|\eta'_1|}{1-|\eta'_1|}.
\end{displaymath}
Since $|\eta'_1|\le 2^{-d'}\le 2^{-1}$ we get
\begin{displaymath}
|\xi|\le 2|\eta'_1|\le 2\times2^{-d'}\le 2^{-d}.
\end{displaymath}
\end{proof}

\subsection{Complex Numbers.}
In an MPFP system a representable complex number $a$ is given by a
pair of representable real numbers $a=x+iy$. \(Here the symbol $+$
is not an operation but a convenient way of representing the pair.\)

We define $\circ_d(a):= \circ_d(x)+i\circ_d(y)$.

\begin{proposition}\label{roundcomplex} Let   $a$ be
a complex number. Then if $a'=\circ_d(a)$, we have $a'=a(1+\eta_1)$,
where $\eta_1$ denotes  a complex number with $|\eta_1|\le
2^{-d}$.
\end{proposition}
\begin{proof}
We may assume that $a\ne0$ ($\,$if $a=0$ then also  its truncation is  $0$
and the claim is trivial$\,$).

Let $a=x+iy$ and $a'=\circ_d(a)=x'+iy'$. Then there exist real numbers
$\varepsilon$ and $\delta$ with $|\varepsilon|\le 2^{-d}$ and
$|\delta|\le 2^{-d}$ such that $x'=x(1+\varepsilon)$ and
$y'=y(1+\delta)$. Then we will have
\begin{displaymath}
a'=x'+iy'=x+iy+\varepsilon x+i\delta y=a\Bigl(1+\frac{\varepsilon
x+i\delta y}{x+iy}\Bigr)=a(1+\eta_1)
\end{displaymath}
and
\begin{displaymath}
|\eta_1|=\Bigl|\frac{\varepsilon x+i\delta
y}{x+iy}\Bigr|=\frac{(\varepsilon^2x^2+\delta^2y^2)^{1/2}}{(x^2+y^2)^{1/2}}\le
\frac{2^{-d}(x^2+y^2)^{1/2}}{(x^2+y^2)^{1/2}}\le 2^{-d}.\qedhere
\end{displaymath}
\end{proof}

In the sequel we will use $\oplus$, $\ominus$,  $\otimes$ and
$\oslash$  to denote the operations performed in the MPFP system.  For
complex numbers these operations are defined in terms of
operations on real numbers in the following way:
\begin{multline*}
(x+iy)\oplus(u+iv):=x\oplus u + i y\oplus v,\quad (x+iy)\ominus(u+iv):=x\ominus u + i y\ominus v,\\
(x+iy)\otimes(u+iv):=(x\otimes u)\ominus(y\otimes
v)+i\bigl\{(x\otimes v)\oplus(y\otimes u)\bigr\},\\
(x+iy)\oslash(u+iv)=((x+iy)\otimes(u-iv))\oslash((u\otimes u)\oplus(v\otimes v)).
\end{multline*}

\begin{proposition}
Let $a$ and $b$ be rounded complex numbers. Then  there exist complex
numbers $\eta_1$ such that
\begin{align}
&a\oplus b=(a+b)(1+\eta_1), \quad & |\eta_1|\le 2^{-d}\label{firstequation}\\
&a\ominus b=(a-b)(1+\eta_1),\quad & \quad |\eta_1|\le 2^{-d}.
\end{align}
\end{proposition}

\begin{proof}
Since the two assertions are similar we only prove the first one.
If $a+b=0$ there is nothing to prove, since in this case  $a\oplus b=0$.

Let $a=x+iy$ and $b=u+iv$. Since we assume that $a$ and $b$ are rounded, there exist $\varepsilon$ and
$\delta\in\R$ with $|\varepsilon|$, $|\delta|\le 2^{-d}$ and such
that
\begin{multline*}
a\oplus b=x\oplus u + i y\oplus v=(x+u)(1+\varepsilon)+i(y+v)(1+\delta)=\\ =
(a+b)\Bigl(1+\frac{\varepsilon(x+u)+i\delta(y+v)}{(x+u)+i(u+v)}\Bigr)=(a+b)(1+\eta)
\end{multline*}
and
\begin{displaymath}
|\eta|^2=\frac{\varepsilon^2(x+u)^2+\delta^2(y+v)^2}{(x+y)^2+(y+v)^2}\le
2^{-2d}.\qedhere
\end{displaymath}
\end{proof}

\begin{proposition}
Let $x$ be a real and $a$ a complex  number, both rounded. Then we have
\begin{align}
x\otimes a & = xa(1+\eta_1), & a\otimes x &= a x(1+\eta_1).
\end{align}
If $a$ and $b$ are rounded complex numbers, then
\begin{equation}
a\otimes b=ab(1+\eta_3).
\end{equation}
\end{proposition}

\begin{proof}
To prove the first  assertion we
may assume that $a\ne0$.  By Proposition \ref{IEEEstandard}
we have numbers $\varepsilon$ and $\delta$ of absolute value less than or
equal to $2^{-d}$ such that
\begin{multline*}
x\otimes a=x\otimes(u+iv)=x\otimes u + i x\otimes v=xu(1+\varepsilon)+ixv(1+\delta)=\\
=xa\Bigl(1+\frac{\varepsilon u+i\delta
v}{(u+iv)}\Bigr)=xa(1+\eta_1)
\end{multline*}
and
\begin{displaymath}
|\eta_1|^2=\frac{\epsilon^2 u^2+\delta^2v^2}{u^2+v^2}\le 2^{-2d}.
\end{displaymath}

Now we consider the case of two complex numbers $a=x+iy$ and
$b=u+iv$.  As usual we may assume that $ab\ne0$.  By the definition
 of $a\otimes b$ we will have
\begin{displaymath}
a\otimes b =\bigl\{x\otimes (u+iv)\}\oplus\bigl\{y\otimes(-v+iu)\}
\end{displaymath}
so that by the above results there exist complex numbers $\delta$, $\xi$
and $\eta$ of absolute value less than $2^{-d}$ such that
\begin{multline*}
a\otimes b=(x(u+iv)(1+\delta))\oplus(y(-v+iu)(1+\xi))=\\ =
\{x(u+iv)(1+\delta)+y(-v+iu)(1+\xi)\}(1+\eta)=\\
ab\Bigl(1+\frac{\delta x(u+iv)+\xi y(-v+iu)}{ab}\Bigr)(1+\eta)=
ab\Bigl(1+\frac{\delta x+ i\xi y}{a}\Bigr)(1+\eta)=\\
=ab(1+\lambda)(1+\eta)
\end{multline*}
where by the Schwarz inequality
\begin{displaymath}
|\lambda|\le \frac{(|\delta|^2+|\xi|^2)^{1/2}|a|}{|a|}\le \sqrt{2}
\; 2^{-d}
\end{displaymath}
Now define $\mu$ such that $1+\mu=(1+\lambda)^{1/2}$. Since $|\lambda|\le \sqrt{2} \,2^{-d}\le 2^{-1/2}$, we can take
\begin{displaymath}
\mu=\sum_{n=1}^\infty \binom{1/2}{n}
\lambda^n=\lambda\sum_{n=1}^\infty \binom{1/2}{n} \lambda^{n-1}
\end{displaymath}
so that
\begin{displaymath}
|\mu|\le \sqrt{2}\,
2^{-d}\sum_{n=1}^\infty\Bigl|\binom{1/2}{n}\Bigr| 2^{-(n-1)/2}\le
0.917608\, 2^{-d}< 2^{-d}.
\end{displaymath}
It follows that
\begin{displaymath}
(1+\lambda)(1+\eta)=(1+\mu)^2(1+\eta)=(1+\eta_3)
\end{displaymath}
since $|\mu|$ and $|\eta|\le 2^{-d}$.
\end{proof}

\subsection{Problem with the definition of the product.}\label{S:2.6}

So, for rounded real numbers $x$ and $y$ we have $x\otimes
y=(xy)(1+\eta_1)$, but for complex numbers we only have  $a\otimes
b=(ab)(1+\eta_3)$.  This  problem is due to the definition  of the
product:
\begin{displaymath}
(x+iy)\otimes(u+iv):=(x\otimes u)\ominus(y\otimes
v)+i\bigl\{(x\otimes v)\oplus(y\otimes u)\bigr\}.
\end{displaymath}

In particular this does not follows the IEEE standard.  The result
of a multiplication must be as if computed exactly and then rounded.
In this vein we  we would define  $a\otimes b=\circ(ab)$. Then by
Proposition \ref{roundcomplex} we would get
\begin{displaymath}
a\otimes b = ab(1+\eta_1),  \qquad \text{($\,$not true, but desirable$\,$).}
\end{displaymath}

It would be desirable to implement the product of complex numbers
satisfying this desideratum.

\subsection{Turing notation.}
Let $A$ be a complex number and $\varepsilon>0$.  We denote by $A+\Turing(\varepsilon)$ a representable number $a$ such that $|a-A|<\varepsilon$. In pseudocode we will write
\begin{equation}
\texttt{a = A + Turing(eps)}
\end{equation}
to indicate that we have applied an algorithm to compute $A$ with an error less than \texttt{eps}.
In this paper we frequently deal with the problem of computing
$A+\Turing(\varepsilon)$ where  $A$  is a given real or complex number.

\subsection{Well implemented functions.}

We will say that a function $f(x)$ is well implemented  in a
domain $\Omega$ if for each representable number $x\in\Omega$ and
with a working precision equal to $d$  the code \texttt{u = f(x)}
gives us a representable number such that $|f(x)-u|\le 2^{-d}$.

In \text{mpmath} we  have the well implemented functions \texttt{log},
\texttt{exp}, \texttt{sqrt}, \texttt{sin}, \texttt{cos} in the
ranges where we have to make use of them.

\subsection{Simple Bounds.}

Sometimes we need a simple bound of a representable number. In
\texttt{mpmath} we can get it easily with the function \texttt{mag}.
For $a\ne0$ real the code \texttt{b=mag(a)} gives us an integer $b$ such that
$2^{b-1}\le |a| <2^b$. For a complex $a$  we will have $|a| <2^b$, it is not 
guaranteed that $b$ is an optimal bound, but it will never be too large by more
than $2$.

\section{The Riemann-Siegel procedure.}

As explained in \cite{A86} the computation of $\zeta(s)$ or $Z(t)$ 
is reduced to the computation of the integral $\Rzeta(s)$ with $t:=\Im  s>0$\footnote{For notations
not explained here, see \cite{A86}.}.
So we want to compute $\Rzeta(s)+\Turing(\varepsilon)$. That is we want
to determine a dyadic number that approximates $\Rzeta(s)$ with an
error less than $\varepsilon$.

The data are: three dyadic numbers $\sigma=\Re s$,  $t=\Im s$ and
$\varepsilon$. We assume that $t$ is sufficiently large and positive.
We also assume that $0<\varepsilon<1$.  (TO DO Complete this).

In practice we will always take $\varepsilon= 2^{-d}$ where $d$ is the working precision \texttt{wpinitial} at the moment we ask for the value of
$\Rzeta(s)$.

The program starts  copying the value of the working
precision in the variable
\texttt{wpinitial = mp.prec}.

With these data we obtain the numbers defined in \cite{A86}*{equation (5)}:  $a=\sqrt{t/2\pi}$, $N=\lfloor a\rfloor$ and
$p=\bigl\{1-2(a-N)\bigr\}$. (Observe that $N$ is \emph{numerically}
well defined since $t$ is a dyadic number.) We also compute
$a^{\sigma}$.

\subsection{End of Computation.}

By the Riemann-Siegel formula
\begin{displaymath}
\Rzeta(s)=\sum_{n=1}^N\frac{1}{n^s}+(-1)^{N-1}U a^{-\sigma}\Bigl(
\sum_{k=0}^K \frac{C_k(q)}{a^k}+RS_K\Bigr)
\end{displaymath}
where $N$ and $K$ are adequately chosen integers and $RS_K$ is the error term.

Roughly the procedure consist in   
determining a simple bound $A_1$  such that $a^\sigma\ge A_1>0$, and then
compute the two numbers
\begin{equation}\label{defS1S2}
S_1:=\sum_{n=1}^N\frac{1}{n^s}+\Turing(\varepsilon/6),\qquad
S_2:=\sum_{k=0}^K \frac{C_k(q)}{a^k}+RS_K+\Turing(A_1\varepsilon/3).
\end{equation}
Once computed these numbers we obtain a simple bound $A_2> |S_2|$ and then compute
\begin{equation}\label{defS3}
S_3:=(-1)^{N-1}U a^{-\sigma}+\Turing(\varepsilon/3A_2).
\end{equation}
Finally we end the computation with the code:
\bigskip
\bigskip
{\footnotesize
\begin{Verbatim}[numbers= left,numbersep=2pt,numberblanklines=false,frame=single,
label=\fbox{\Large End of computation}]

(* we have computed  S1, S2, S3  *)

mp.prec = 15
absS1 = abs(S1); absS2 = abs(S2*S3)
d = max( 6, d0 + mag( 6 * (3 * absS1 + 7 * absS2) ) )
mp.prec = d
S1' = round(S1); S2' = round(S2); S3' = round(S3)
R = S1' + ( S3' * S2' )
return R
\end{Verbatim}
}

In the first four lines of this code  we compute the precision $d$ at which we
make the next computation. Of course the computation of this precision is done 
only to  15 digits of precision. This explain the second line.
Lines 6 and 7 are not real.  In the true code we will only write 
\texttt{R = S1 + ( S3 * S2 )}, but according with the IEEE standard this 
is equivalent to first a rounding (line 6), the computation and then a
final rounding.

We have to show that $R=\Rzeta(s)+\Turing(\varepsilon)$. First observe
that by \eqref{defS1S2}, \eqref{defS3}, the fact that $|U|\le 1$, and
the Riemann Siegel formula \cite{A86}*{(4)}  we have
\begin{multline*}
S_1+S_3S_2=\sum_{n=1}^N\frac{1}{n^s}+\Turing(\varepsilon/6)+\\+
\bigl\{(-1)^{N-1}U a^{-\sigma}+\Turing(\varepsilon/3A_2)\bigr\}\{
\sum_{k=0}^K \frac{C_k(q)}{a^k}+RS_K+\Turing(A_1 \varepsilon/3)\bigr\}=\\
=\Rzeta(s)+\Turing(\varepsilon/6)+\Turing(\varepsilon/3)+\Turing(\varepsilon/3)=
\Rzeta(s)+\Turing(5\varepsilon/6).
\end{multline*}

Since $d\ge6$ we have $6\cdot 2^{-d}<0.1$ so that   $|\eta_r|\le
1.06\cdot r \cdot 2^{-d}$ for $1\le r\le 6$  by Lemma \ref{boundExp}. Therefore
\begin{multline*}
|S'_1\oplus(S'_3\otimes S'_2)-\Rzeta(s)|\le \frac{5\varepsilon}{6}+
|S'_1\oplus(S'_3\otimes S'_2)-(S_1+S_2S_3)|=\\
=\frac{5\varepsilon}{6}+|S'_1\oplus(S'_3  S'_2)(1+\eta_3)-(S_1+S_2S_3)|=\\
=\frac{5\varepsilon}{6}+|(S'_1+(S'_3 S'_2)(1+\eta_3))(1+\eta_1)-(S_1+S_2S_3)|=\\
=\frac{5\varepsilon}{6}+|(S_1(1+\eta_1)+(S_3 S_2)(1+\eta_1)(1+\eta_1)(1+\eta_3))(1+\eta_1)
-(S_1+S_2S_3)|=\\
=\frac{5\varepsilon}{6}+|(S_1(1+\eta_2)+(S_3 S_2)(1+\eta_6)
-(S_1+S_2S_3)|=\\
\frac{5\varepsilon}{6}+|\eta_2 S_1+\eta_6 S_2S_3|\le
\frac{5\varepsilon}{6}+2^{-d} (3|S_1|+7|S_2S_3|)<\varepsilon
\end{multline*}
since the election of $d$ implies that $2^{-d} (3|S_1|+7|S_2S_3|)<\frac{\varepsilon}{6}$.

Define $\varepsilon_1:=\varepsilon/6$ and
$\varepsilon_2:=A_1\varepsilon/3$.

\subsection{Number of terms in the Riemann-Siegel Correction.}

The main remaining problem is to compute the Riemann-Siegel sum
\begin{equation}\label{sumK}
\sum_{k=0}^{K} \frac{C_k(p)}{a^k}+RS_K+\Turing(\varepsilon_2).
\end{equation}
In principle $K$ could be any natural  number. The problem is that
we do not know how to compute the rest $RS_K$; we only have bounds.
So we must choose $K$ in such a way that the absolute value of
$RS_K$ is small.

We have better bounds for the terms of the sum in \eqref{sumK} than
for the rest, so that as in many similar situations,
it is advantageous to choose $0<L<K$ and to put
\begin{displaymath}
\sum_{k=0}^{K} \frac{C_k(p)}{a^k}+RS_K=\sum_{k=0}^{L-1}
\frac{C_k(p)}{a^k}+\sum_{k=L}^{K} \frac{C_k(p)}{a^k}+RS_K
\end{displaymath}
and find a  bound of the new rest
\begin{displaymath}
\sum_{k=L}^{K} \frac{C_k(p)}{a^k}+RS_K.
\end{displaymath}
In this way we will get almost the same result as by applying the usual
rule of   thumb: \emph{the error is of the order of magnitude of the
first term omitted}.

This is the content of the following Theorem.

\begin{theorem}\label{propL}
Let   $L\ge1$ be an integer such that
\begin{equation}\label{maincondition}
3c\frac{\Gamma(L/2)}{(b a)^L}<\varepsilon_2
\end{equation}
where $b=b(\sigma)$ and $c=c(\sigma)$ are defined in \cite{A86}*{Theorem 4.1,  (4.2)}.
Assume also that
\begin{equation}\label{secondarycondition}
3L+4<\frac{8}{25}a^2=\frac{4t}{25\pi}\quad\text{and}\quad
3L+2+\sigma\ge0.
\end{equation}
Then with an adequate choice of K we will have
\begin{displaymath}
\sum_{k=0}^{K}
\frac{C_k(p)}{a^k}+RS_K+\Turing(\varepsilon_2)=\sum_{k=0}^{L-1}\frac{C_k(p)}{a^k}
+\Turing(\varepsilon_2/2)
\end{displaymath}
in the sense that every solution to the problem in the right hand side will be a solution of the problem in the left hand side.
\end{theorem}

\begin{proof}
Assume we have $L$ satisfying the conditions  \eqref{maincondition}
and \eqref{secondarycondition}. Taking $K=3L+4$ we have
\begin{displaymath}
\sum_{k=0}^{K} \frac{C_k(p)}{a^k}+RS_K=\sum_{k=0}^{L-1}
\frac{C_k(p)}{a^k}+\sum_{k=L}^{K} \frac{C_k(p)}{a^k}+RS_K
\end{displaymath}
and by Theorems 4.1  and 4.2 of \cite{A86} we will have ($\,$since
$3L+4+\sigma\ge2$, Theorem 4.2 of \cite{A86} applies$\,$)
\begin{displaymath}
T:=\Bigl|\sum_{k=L}^{K} \frac{C_k(p)}{a^k}+RS_K\Bigr|\le
c\sum_{k=L}^{K}\frac{\Gamma(k/2)}{(ba)^k}+c_1
\frac{\Gamma((K+1)/2)}{(b_1a)^{K+1}}
\end{displaymath}
where $b$, $c$, $b_1$ and $c_1$ are the coefficients appearing in
Theorems 4.1 and 4.2 of \cite{A86} ($\,$recall that they depend on
$\sigma$$\,$).

For $x\ge1$ we have $\Gamma((x+1)/2)\le \Gamma(x/2)\sqrt{x/2}$.
Observe also that $b\ge2$. Therefore the quotient of two
consecutive terms of the first sum is
\begin{displaymath}
\frac{\Gamma((k+1)/2)}{ba\Gamma(k/2)}\le
\frac{\sqrt{k/2}}{ba}\le \frac{\sqrt{K/2}}{ba}=
\frac{\sqrt{(3L+4)/2}}{ba}<\frac{1}{ba}\Bigl(\frac{4}{25}a^2\Bigr)^{1/2}\le\frac{1}{5}.
\end{displaymath}
Then we get
\begin{displaymath}
T\le
c\frac{\Gamma(L/2)}{(ba)^L}\Bigl(1+\frac{1}{5}+\frac{1}{5^2}+\cdots+\frac{1}{5^{K-L}}+
\frac{c_1}{c}\frac{(b/b_1)^{K+1}}{5^{K-L+1}}\Bigr).
\end{displaymath}

Since  $c_1/c\le 0.5 (0.9)^{\lceil
-\sigma\rceil}/(2^{-\sigma}/\pi\sqrt{2})\le \pi/\sqrt{2}\approx
2.22144$ for $\sigma<0$ and
$c_1/c=(2^{3\sigma/2}/7)/(9^\sigma/\pi\sqrt{2})\le
\pi\sqrt{2}/7\approx 0.634698$ for $\sigma>0$ we will have (for
$L\ge1$)
\begin{displaymath}
\frac{c_1}{c}\frac{(b/b_1)^{K+1}}{5^{K-L+1}}\le
2.23\frac{(2.26*1.1)^{K+1}}{5^{K-L+1}}=2.23\frac{(2.26*1.1)^{3L+5}}{5^{2L+5}}<0.04165.
\end{displaymath}
It follows that
\begin{displaymath}
\Bigl|\sum_{k=L}^{K-1} \frac{C_k(p)}{a^k}+RS_K\Bigr|\le
c\frac{\Gamma(L/2)}{(ba)^L}(5/4+ 0.05)<\varepsilon_2/2.\qedhere
\end{displaymath}
\end{proof}

\subsection{What precision may we get?}\label{S:3.3}
Proposition \ref{propL} determines  what precision we may  get with the Riemann-Siegel expansion.  We will speak about the computation of  $\Rzeta(\sigma+it)$ but the considerations extends to the zeta function.  We shall consider $\sigma$ fixed and $t\gg1$.

We will choose  $A_1:=\frac12 2^{\texttt{mag}(a^\sigma)}\le a^\sigma<2A_1$. 
Then $\varepsilon_2:= A_1\varepsilon/3$ satisfies
\begin{equation}
\frac{\varepsilon a^\sigma}{6}<\varepsilon_2\le\frac{\varepsilon a^\sigma}{3}.
\end{equation}
For $r$ not too small, the minimum value  of $\Gamma(x)/r^x$ is approximately taken at the point $x=r$, and is approximately equal to $\sqrt{2\pi/ r}\, e^{-r}$.  Then we will get a value of $L$ satisfying \eqref{maincondition} if
\begin{displaymath}
\sqrt{2\pi/b^2 a^2}e^{-b^2 a^2}<\frac{\varepsilon_2}{3c}
\end{displaymath}
Here $b$ and $c$  are functions of $\sigma$ \(determined in \cite{A86}*{(26)}\).
Since we assume $a$ big, there will be a value of $L$ satisfying the above condition if
\begin{displaymath}
\exp\Bigl(-\frac{b^2}{2\pi}t\Bigr)<\frac{\varepsilon a^\sigma}{18c}
\end{displaymath}

For example in the case of $\sigma=1/2$,  $b=2$, $c=3/\pi\sqrt{2}$, we can assume $a^\sigma>18c$ and the above condition is implied by $e^{-2t/\pi}<\varepsilon$.

Thus we can compute $\zeta(\tfrac12 +it)$ with error less than $\varepsilon$ by the Riemann-Siegel formula if $e^{-2t/\pi}<\varepsilon$. In other words
we can get $d$ binary digits if $\frac{2}{\pi \log 2} t>d$.

Since the Riemann-Siegel method is only useful for $t$ big \(in other case there are better methods\) in practice it is difficult to imagine a situation where the condition $\frac{2}{\pi \log 2} t>d$ will be a real problem.

There is other problem that is more important. If we want to compute $\zeta(s)$
with $d\sim t$ digits, then  the number $L$ is big.  This implies that a large 
number of terms of the Riemann-Siegel corrections must be computed,
 and these terms are difficult to compute. 
In condition  \eqref{secondarycondition}  we see that the error in the 
Riemann-Siegel formula is small even when $L\sim t$. The problem will be that
we have to compute so many terms of the expansion that the advantage over Mac-Laurin
formula will not be clear in these cases.

In fact in Proposition \ref{PropComputPolynomial} we shall need to substitute the first condition in \eqref{secondarycondition} by
\begin{equation}\label{thirdcondition}
3L<\frac{2a^2}{25}=\frac{t}{25\pi}.
\end{equation}
By the above considerations it is clear that this is a very modest restriction.

\subsection{How to get $\sum_{n=1}^N a_n+\Turing(\varepsilon)$. }\label{sumsection}

We may not assume that the law $x\oplus y$ is commutative or  associative. ($\,$In fact the usual operations defined on computers are commutative but not associative ($\,$see Knuth \cite{Knuth}$\,$)$\,$).
But for convenience  we will use the notation
\begin{displaymath}
\bigoplus_{j=1}^n x_j = (x_1\oplus x_2\oplus \cdots \oplus x_n)= (\cdots ((x_1\oplus x_2)\oplus x_3)\oplus \cdots \oplus x_{n-1})\oplus x_n
\end{displaymath}
to denote the computed sum when the additions takes place in the order in which they are written.

The following Proposition can be found in \cite{Wilkinson}, and can be proved easily by induction.

\begin{proposition}\label{sumofrounded}
Let $x_k$ be rounded real or complex numbers, then
\begin{equation}\label{rulesum}
x_1\oplus x_2\oplus \cdots \oplus x_n
=x_1(1+\eta_n)+x_2(1+\eta_{n-1})+\cdots +x_n(1+\eta_1).
\end{equation}
\end{proposition}

Note that the error obtained is dependent on the order of
summation. The upper bound for the error is smallest if the terms
are added in order of increasing absolute magnitude  since then the
largest factor $(1+\eta_n)$ is associated with the smallest $x_i$.

Assume that we want to compute
\begin{equation}
S:=\sum_{n=1}^N a_n +\Turing(\varepsilon).
\end{equation}
We will assume that for each $\delta>0$ and $n\le N$ we know how to get
$a_n+\Turing(\delta)$, and that we know a bound $A$ such that $|a_n|\le A$ for each $n\le N$.
\medskip

Then we will use the following procedure
\bigskip
\bigskip
{\footnotesize
\begin{Verbatim}[numbers= left,numbersep=2pt,numberblanklines=false,frame=single,
label=\fbox{\Large Program to compute $\sum_{n=1}^N a_n+\Turing(\varepsilon)$}]

if N * A < eps:
    return 0

delta = eps/(2*N)
v={}
for n in range(1,N+1):
    v[n]=a[n]+Turing(delta)

mp.prec = 15
d= max( mag(10*(N+1)), mag(2.2*A*(N+3)**2 / eps) ) + 1
mp.prec = d

sum = 0
for n in range(1, N+1):
    sum = sum + v[n]
return sum
\end{Verbatim}
}

\begin{proof}[Proof of the correctness of the algorithm]
It is clear that if $NA<\varepsilon$, then $0$ is a representable number and
$\left|0-\sum_{n=1}^N a_n\right|<\varepsilon$. This explains line 1--2.
If we pass this line then $NA\ge \varepsilon$, therefore in what follows $A>\delta$.

In lines 3--6 we get representable numbers $v_n$ such that $|v_n-a_n|<\varepsilon/2N$, and by assumption we know how to find these numbers. Then we can put $v_n=a_n+\alpha_n$ with
$|\alpha_n|\le \delta$

In line 7--9 we fix the working precision to a number $d$ such that
\begin{displaymath}
10(N+1)< 2^d,\qquad \frac{2.2\; A(N+3)^2}{\varepsilon} < 2^d.
\end{displaymath}
(The $+1$ in line 8 guarantees that  we get these inequalities even when
we take account of the factor $(1+\eta_4)\le2$  that appears when we
perform the computation \texttt{2.2*A*(N+3)**2 / eps}.)

In line 10--13 observe that $\texttt{sum}$ is always rounded, so by Proposition \ref{sumofrounded} we will get at the end of the for loop
\begin{align*}
\texttt{sum}&=v'_1(1+\eta_N)+v'_2(1+\eta_{N-1})+\cdots +v'_n(1+\eta_1)\\
&=v_1(1+\eta_{N+1})+v_2(1+\eta_{N})+\cdots +v_N(1+\eta_2)\\
&=(a_1+\alpha_1)(1+\eta_{N+1})+(a_2+\alpha_2)(1+\eta_{N})+\cdots +(a_N+\alpha_N)(1+\eta_2)
\end{align*}
so that, by Lemma \ref{boundExp},
\begin{align*}
\Bigl|\texttt{sum}-\sum_{n=1}a_n\Bigr|&\le \sum_{n=1}^N |a_n\eta_{N-n+2}|+\delta\sum_{n=2}^{N+1}|1+\eta_n|\le
\\ &\le (A+\delta)*1.06\cdot \frac{N(3+N)}{2}2^{-d}+N\delta<\\
&<\frac{2.2 A(N+3)^2 2^{-d}}{2}+\frac{\varepsilon}{2}<\varepsilon.\qedhere
\end{align*}

\end{proof}

To finish we mention again that when computing a sum we should try to sum the terms
in increasing order of their absolute values.

\subsection{Computing the Riemann-Siegel sum.}
We will compute the sum
\begin{displaymath}
\texttt{rssum} := \sum_{k=0}^{L-1}\frac{C_k(p)}{a^k}+\Turing(\varepsilon_2/2).
\end{displaymath}
In Theorem \ref{propL} we have seen that  the terms of this sum are mainly decreasing.
So we will compute the sum in the form
\begin{equation}
\texttt{rssum} = \sum_{k=1}^{L}\frac{C_{L-k}(p)}{a^{L-k}}+\Turing(\varepsilon_2/2).
\end{equation}
By Theorem \cite{A86}*{Theorem 2, (29)} and Theorem \ref{propL} we know that the 
terms of this sum are bounded by $c\Gamma(k/2)/(ba)^k$ for $k\ge1$ so that 
they are bounded by
\begin{equation}
\texttt{rsbound}:=c\sqrt{\pi}/b a.
\end{equation}
So, applying the procedure of Section \ref{sumsection} we will have to solve the problem
\begin{equation}
\texttt{term[k]} := \frac{C_k(p)}{a^k}+\Turing(\varepsilon_3)\quad\text{where}\quad \varepsilon_3=\frac{\varepsilon_2}{4L}.
\end{equation}
and perform the summation with  precision
\begin{equation}\label{defwprssum}
\texttt{wprssum}=\max\Bigl\{\texttt{mag}(10(L+1)),\texttt{mag}\Bigl(\frac{4.4c\sqrt{\pi}}
{\varepsilon_2 b a}(L+3)^2\Bigr)\Bigr\}+1.
\end{equation}

\subsection{Computing \texttt{term[k]}. }
By \cite{A86}*{equation (39)}
\begin{equation}\label{oneterm}
 \frac{C_k(p)}{a^k}=\frac{1}{(\pi^2 a)^{k}}\sum_{j=0}^{\lfloor3k/2\rfloor}
\Bigl(\frac{\pi}{2i}\Bigr)^j d^{(k)}_j
F^{(3k-2j)}(p)
\end{equation}

In the proof of Proposition 6 in \cite{A86} we obtained bounds for the terms of the sum
in \eqref{oneterm}
that were increasing. So, we will sum in the indicated order. By computing the maximum of these bounds we obtained in \cite{A86}*{Proposition 6} a bound $T_k/a^k$ for all the terms of this sum.

With this bound we can apply the procedure given in Section \ref{sumsection}.  The sum has
$1+\lfloor 3k/2\rfloor \le 1+3k/2$ terms. So we must compute each term with an error less than
$\varepsilon_3/(3k+2)$. We prefer to compute a little more precise so that our epsilon does not depend on $k$.  Since $k\le L-1$ we have $3k+2\le 3L$.  So, we must solve the problem
\begin{equation}\label{tkjplusTuring}
\frac{t^{(k)}_j}{a^k}+\Turing(\varepsilon_4)\quad\text{where}\quad
t^{(k)}_j :=
\Bigl(\frac{\pi}{2i}\Bigr)^j \frac{d^{(k)}_j}{\pi^{2k}}
F^{(3k-2j)}(p), \quad\varepsilon_4=\frac{\varepsilon_3}
{3L}.
\end{equation}

By the general result of Section \ref{sumsection}  must compute the sum in \eqref{oneterm} with a working precision equal to the least natural number $d$ such that
\begin{displaymath}
2^d >10(3k/2+1),\qquad 2^d> 2.2\frac{T_k}{a^k}(3k/2+4)^2\frac{1}{\varepsilon_3}.
\end{displaymath}
We have
\begin{displaymath}
\frac{3k}{2}+2<\frac{3k}{2}+4\le \frac{3L-3}{2}+4<\frac{3(L+2)}{2}.
\end{displaymath}
So the first condition on $d$ is satisfied if we take $d>\texttt{mag}(40(L+2))$.

For the second condition we have
\[
2.2\frac{T_k}{a^k}(3k/2+4)^2\frac{1}{\varepsilon_3}<2.2\frac{9}{4}
\frac{(L+2)^2\,A}{\varepsilon_3}\frac{\Gamma(k+\frac12)^{1/2}}{(B_1 a\sqrt{\pi})^k}
<68\frac{(L+2)^2\, A}{\varepsilon_3}\frac{\Gamma(k+\frac12)^{1/2}}{(B_1 a\sqrt{\pi})^k}
\]
So that we can take the precision as \texttt{wpterm[k]}
\begin{equation}\label{defwpterm}
\texttt{wpterm[k]} = \max\Bigl\{\texttt{mag}(40(L+2)),
\texttt{mag}\Bigl(68\frac{(L+2)^2\, A}{\varepsilon_3}\frac{\Gamma(k+\frac12)^{1/2}}{(B_1 a\sqrt{\pi})^k}\Bigr)\Bigr\}.
\end{equation}
The reason for the choice of the two constant $40$ and $68$ will appear in Section \ref{sectiontcoef}. We shall compute wpterm[k] at the same time that we compute wptcoef[k] and these choices simplify the simultaneous computation.

\subsection{Computing $a(1+\eta_1)$.}\label{trunca}

Recall that $a=\sqrt{t/2\pi}$. Given the value $d$  of \texttt{mp.prec}
we want to compute an approximate value $a(1+\eta_1)$.
This is achieved by the following procedure:
\bigskip
\bigskip
{\footnotesize
\begin{Verbatim}[numbers= left,numbersep=2pt,numberblanklines=false,frame=single,
label=\fbox{\Large Computing $a(1+\eta_1)$.}]

def trunc_a(t):
    wp = mp.prec
    mp.prec=wp+2
    aa = sqrt(t/(2*pi))
    mp.prec=wp
    return(aa)
    \end{Verbatim}
}

The value of the variable $\texttt{aa}$ would be
\begin{displaymath}
\sqrt{\frac{t(1+\eta_1)}{2\pi(1+\eta_1)}(1+\eta_1)}(1+\eta_1)=a\sqrt{\frac{1+\eta_2}{1+\eta_1}}
(1+\eta_1).
\end{displaymath}
\(Recall that each instance of $\eta_1$ may represent a different number.\)
Let $d$ be the initial value of \texttt{BynaryPrec} and $d'=d+2$ the actual value.
Then we will have
\begin{align*}
|\sqrt{1+\eta_1}-1|&= \Bigl|\eta_1\sum_{n=1}^\infty \binom{1/2}{n}\eta_1^{(n-1)}\Bigr|\le
0.585787\times 2^{-d'}\\
\Bigl|\frac{1}{\sqrt{1+\eta_1}}-1\Bigr|&=
\Bigl|\eta_1\sum_{n=1}^\infty \binom{-1/2}{n}\eta_1^{(n-1)}\Bigr|\le 0.828428\times 2^{-d'}
\end{align*}
so that $\texttt{aa}= a(1+\varepsilon)$ with
\begin{displaymath}
|\varepsilon|< (1+\frac{0.586}{4}\times 2^{-d})^2(1+\frac{0.829}{4}\times 2^{-d})(1+\frac{1}{4}2^{-d})-1\le 0.984\times 2^{-d}<2^{-d}
\end{displaymath}
\(the factor $0.984$ appears when we expand the products and  substitute all the powers of $2^{-d}$
 by $2^{-d}$\).

\subsection{Computing the quotients $t^{(k)}_j/a^k+\Turing(\varepsilon_4)$.}

We have to compute the powers of $a^{-1}$. Then we multiply $t^{(k)}_j \times (a^{-1})^k$.
In \cite{A86}*{Proposition 6} we have given a bound $|t^{(k)}_j|\le T_k$.
We shall assume that we have computed
\begin{equation}\label{pbtcoef}
\texttt{tcoef[k,j]}=t^{(k)}_j+\Turing(\delta_k).
\end{equation}
Our problem is to determine to which precision we have to compute the powers $a^{-k}$,
how to choose
the numbers $\delta_k$ and which working precision to use in the
computation of the products $t^{(k)}_j \times (a^{-1})^k$.

\begin{proposition}
In order to compute $\texttt{tv[k,j]}=t^{(k)}_j+\Turing(\varepsilon_4)$ for $0\le k<L$, $0\le j\le \lfloor 3k/2\rfloor$ we determine for each $k$ a working precision \texttt{wptv[k]} as
the least natural number $d$ such that
\begin{equation}
2^d>10(2k+3),\qquad 2^d>68\frac{L A(\sigma)}{\varepsilon_4
}\frac{\Gamma(k+1/2)^{1/2}}{(B_1 a \sqrt{\pi})^k}.
\end{equation}
Also we compute $\texttt{tcoef[k,j]}= t^{(k)}_j+\Turing(\varepsilon_4 a^k/4)$ and
then follows the following procedure:
\end{proposition}

\bigskip
\bigskip
{\footnotesize
\begin{Verbatim}[numbers= left,numbersep=2pt,numberblanklines=false,frame=single,
label=\fbox{\Large Computing the quotients $t^{(k)}_j/a^k+\Turing(\varepsilon_4)$.}]

# Computing the powers av[k] = a**(-k)
mp.prec = wptv[0]+2
a = trunc_a(t)  # get an approximate value of a. See Section 3.6
av = {}
av[0] = 1
av[1] = av[0]/a

mp.prec = wptv[0]
for k in range(2,L):
    av[k] = av[k-1] * av[1]

# Computing the quotients
tv = {}
for k in range(0,L):
    mp.prec = wptv[k]
    for ell in range(0,3*k/2 + 1):
        tv[k,ell] = tcoef[k,ell]* av[k]
\end{Verbatim}
}

Recall that $a=\sqrt{a/2\pi}$. We want to compute it only once. Hence we assume that we have
computed it before with a greater precision than is needed.
The variable $\texttt{a}$ represents this approximate value.
So, $\texttt{a}=a(1+\eta_1)$, where $|\eta_1|\le 2^{-d-2}$.

From lines 4--5 we get $\texttt{av[1]}=\frac{1}{a(1+\eta_1)}(1+\eta_1)=
a^{-1}(1+2\eta_1)(1+\eta_1)$. When, in line 6 we change the value of
\texttt{mp.prec} to $wptv[0]$, this inequality will become
$\texttt{av[1]}=a^{-1}(1+\eta_1)$. This is true since
$\frac{1+x/4}{1-x/4}<1+x$ for $0<x<1/2$.

Then, in lines 6--8, we get $\texttt{av[k]}=a^{-k}(1+\eta_{2k})$.

In line 11 we change the working precision. We have  $\texttt{wptv[0]}\ge \texttt{wptv[k]}$.
So the new $\eta_k$ are larger. Hence, in line 13 we get
 \(observe that the new working precision by truncation introduces
 two factors $(1+\eta_1)$\)
\begin{displaymath}
\texttt{tv[k,j]}=(t^{(k)}_j+\alpha^{(k)}_j)a^{-k}(1+\eta_{2k+3})\quad\text{where} \quad |\alpha^{(k)}_j|\le \delta_k.
\end{displaymath}
We get the desired result \eqref{tkjplusTuring} if
\begin{displaymath}
\frac{\delta_k}{a^k}(1+\eta_{2k+3})+\frac{T_k}{a^k}\eta_{2k+3}<
\varepsilon_4.
\end{displaymath}
We take $(2k+3)\cdot 2^{-d_k}<0.1$ so that Lemma \ref{boundExp} applies.
It follows that $|\eta_{2k+3}|<1.06\cdot(2k+3) \cdot 2^{-d_k}<1$.
Then the above condition is satisfied if we take
\begin{equation}\label{aprecision}
\delta_k= \frac{\varepsilon_4 a^k}{4} \quad\text{and}\quad
\frac{T_k}{a^k}\cdot 1.06 \cdot(2k+3)\cdot
2^{-d_k}<\frac{\varepsilon_4}{2}.
\end{equation}
By the definition of $T_k$ in \cite{A86}*{(58)} and
\eqref{secondarycondition} we have
\begin{displaymath}
\frac{\frac{T_{k+1}}{a^{k+1}}(2k+5)}{\frac{T_{k}}{a^{k}}(2k+3)}=\frac{(k+1/2)^{1/2}}{B_2a}
\frac{2k+5}{2k+3}\le
\frac{5}{3}\frac{L^{1/2}}{B_2a}\le\frac{5}{3}(8/75)^{1/2}B_2^{-1}<1.
\end{displaymath}
It follows that the \texttt{wptv[k]} are decreasing with increasing $k$.

The second condition on $d_k$ in \eqref{aprecision} is equivalent to
\begin{displaymath}
2^{d_k}>\frac{2.12(2k+3)}{a^k\varepsilon_4}A(\sigma)\frac{\Gamma(k+\frac12)^{1/2}}{ B_2^k}
\end{displaymath}
We have $2k+3\le 2L+1<2L+2$ and $B_2=B_1\sqrt{\pi}$, hence we can take
\begin{displaymath}
\frac{2.12(2k+3)}{a^k\varepsilon_4}A(\sigma)\frac{\Gamma(k+\frac12)^{1/2}}{ B_2^k}\le \frac{4.24 (L +2) A(\sigma)}{\varepsilon_4}\frac{\Gamma(k+\frac12)^{1/2}}{
(aB_1\sqrt{\pi})^k}\le 2^{d_k}
\end{displaymath}
If we compare this with the definition of \texttt{wptcoef[k]}  in \eqref{defwptcoef}, we see that we can take $d_k$ equal to \texttt{wptcoef[k]}.

Summing up we can take $\texttt{wptv[k]}=\texttt{wptcoef[k]}$. That is we compute with three or four more binary digits, saving to compute a new precision for each $k$.

\subsection{Computing $\texttt{tcoef[k,j]}=t^{(k)}_j+\Turing(\delta_k)$. }\label{sectiontcoef}
By \eqref{tkjplusTuring} we see that we can reduce the problem to the computation of the
coefficients $d^{(k)}_j$,  the derivatives $F^{(m)}(p)$ for $0\le m\le 3L-3$ and the
powers  $\pi^{-r}$ for $0\le r\le 2L-2$.  This is the objective of this Section.

In \cite{A86}*{(49), (55)} we have given the bounds $|d^{(k)}_j|\le D^{(k)}_j$ and $|F^{(m)}(p)|\le F_m$.
\smallskip

{\footnotesize
We give the reasons of our choice of epsilons here.
\(What is not strictly needed to the proof of the correctness
of the algorithm will be set in small size.\)

In view of \eqref{tkjplusTuring} and the value of $\delta_k$ in
\eqref{aprecision} a naive application of the rules of the product will require
the computation of
\begin{displaymath}
d^{(k)}_j+\Turing\Bigl(\frac{a^k \varepsilon_4}{2^{-j}\pi^{j-2k}
F_{3k-2j}}\Bigr)\quad\text{and}\quad
F^{(3k-2j)}(p)+\Turing\Bigl(\frac{a^k\varepsilon_4}{2^{-j}\pi^{j-2k}
D^{(k)}_j}\Bigr).
\end{displaymath}
The second of these conditions  is not satisfactory. We want to compute
$F^{(m)}(p)$ only once. So we search the minimum of
\begin{displaymath}
u_{k,j}:=\frac{a^k\varepsilon_4}{2^{-j}\pi^{j-2k} D^{(k)}_j}=
\frac{a^k\varepsilon_4 B_1^k\{(3k-2j)!\}^{1/2}}{2^{-j}\pi^{j-2k} A(\sigma)2^j
\Gamma(k+1/2)^{1/2}}
\end{displaymath}
when $m=3k-2j$.  To get this minimum we compute
\begin{multline*}
\frac{u_{k,j}}{u_{k+2,j+3}}=\frac{1}{\pi a^2 B_1^2}\sqrt{(k+3/2)(k+1/2)}\le\\
 \le\frac{1}{\pi B_1^2
a^2}\sqrt{(L+1/2)(L-1/2)}\le \frac{ L}{\pi B_1^2 a^2}\le \frac{1}{\pi a^2 B_1^2}
\frac{8}{75} a^2<1
\end{multline*}
by \eqref{secondarycondition}.
\(Observe that
$m=3k-2j=3(k+2)-2(j+3)$\).

Therefore, the minimum is obtained when $m=3k-2j$ and $k$ is minimal. Hence for
$k=\lfloor m/3\rfloor+(m\bmod3)\ge m/3$ and $j=m\bmod3$.
It follows that
\begin{multline}\label{ukj}
u_{k,j}\ge \frac{\varepsilon_4(\pi^2A B_1)^{m/3}}{A(\sigma)}\frac{
\{m!\}^{1/2}} {\pi^{m\bmod3} \Gamma(\lfloor m/3\rfloor+(m\bmod3)+1/2)^{1/2}}\ge\\
\ge\frac{\varepsilon_4(\pi^2A B_1)^{m/3}}{\pi^2 A(\sigma)}\Bigl(\frac{m!}{\Gamma(m/3+2)}\Bigr)^{1/2},\qquad
m=3k-2j.
\end{multline}
} 
\medskip

We define
\begin{equation}\label{defeps5}
\widetilde{\varepsilon_5}(m):=\min\Bigl\{4 F_m, \frac{(\pi^2 B_1
a)^{m/3} }{316\,A(\sigma)}\Bigl(\frac{m!}
{\Gamma(m/3+2)}\Bigr)^{1/2}\,\varepsilon_4\Bigr\}.
\end{equation}
\(We want to compute
$F^{(m)}(p)+\Turing(\widetilde{\varepsilon_5}(m))$, since
$|F^{(m)}(p)|\le F_m$ it makes no sense that
$\widetilde{\varepsilon_5}(m)>4 F_m$.\)
\medskip

{\footnotesize
Since we want to compute $2^{-j}\pi^{j-2k}d^{(k)}_jF^{(3k-2j)}(p)$ with an
error less than $\delta_k$ \(see  \eqref{pbtcoef} and \eqref{aprecision}\) roughly  we need to compute with
$d_{k,j}$ binary  digits, where
\begin{displaymath}
2^{-j}\pi^{j-2k}D^{(k)}_jF_{3k-2j} 2^{-d_{k,j}}=
\delta_k=\frac{\varepsilon_4 a^k}{4}.
\end{displaymath}
Since
\begin{displaymath}
\frac{2^{d_{k,j}}}{2^{d_{k,j+1}}}=\sqrt{\frac{3k-2j-1}{3k-2j}}<1
\end{displaymath}
we have $d_{k,j}\le d_{k,\lfloor 3k/2\rfloor}$.  Some examples
indicate that the needed precision $d=d_{kj}$,  for a fixed $k$,
increases smoothly with $j$, and decreases with $k$, so that it
appears reasonable to substitute $d_{k,j}$ for the maximum $d_k$.

In this way we get
\begin{equation}\label{fordk}
2^{d_{k,j}}=\frac{4}{\varepsilon_4 a^k} 2^{-j}\pi^{j-2k}D^{(k)}_jF_{3k-2j}\le
\frac{4 A(\sigma)}{\varepsilon_4
}\frac{\Gamma(k+1/2)^{1/2}}{(B_1 a \sqrt{\pi})^k},\qquad
0\le j\le 3k/2.
\end{equation}
} 

In fact, we will  take as working precision
\begin{equation}\label{defwptcoef}
\texttt{wptcoef[k]}:=\max\Bigl\{\texttt{mag}\Bigl(68\frac{(L +2) A(\sigma)}{\varepsilon_4
}\frac{\Gamma(k+1/2)^{1/2}}{(B_1 a \sqrt{\pi})^k}\Bigr), \texttt{mag}(40(L+2))\Bigr\}.
\end{equation}

Now assume that we have solved the two problems:
\begin{align}
\texttt{Fp[m]}&= F^{(m)}(p)+\Turing(\widetilde{\varepsilon_5}(m)),
\qquad 0\le m\le 3L-3,\label{problemFp}\\
\texttt{d[k,j]}&=d^{(k)}_j+\Turing(\gamma_{k,j}), \qquad 0\le k\le
L-1,\quad 0\le j\le 3k/2,\label{problemd}
\end{align}
where $\widetilde{\varepsilon_5}(m)$ is defined on \eqref{defeps5}
and $\gamma_{k,j}$ will be given in \eqref{defgammak}.

When we compute these coefficients, the powers of $\pi$ needed are
computed to precision $\texttt{wppi} >\texttt{wptcoef[k]}$. They are
stored in the variables \texttt{pipower[k]}. So that we have
$\texttt{pipower[k]}=\pi^k (1+\eta_{2k})$ with the $\eta$'s corresponding to
the precision \texttt{wppi}.

Then the following procedure gives us the \texttt{tcoef[k,j]}.

\bigskip
\bigskip
{\footnotesize
\begin{Verbatim}[numbers= left,numbersep=2pt,numberblanklines=false,frame=single,
label=\fbox{\Large Computing the coefficients \texttt{tcoef[k,ell]}.}]

# computing  the needed wp
# we compute simultaneously wptcoef[k] and wpterm[k]
    wptcoef={}
    mp.prec = 15
    c1 = mag(40*(L+2))
    c2 = mag(68*(L+2)*A)
    c4 = mag(B1*a*math.sqrt(pi))-1
    for k in range(0,L):
        c3 = c2 - k*c4+mag(fac(k+1/2.))/2.
        wptcoef[k] = max(c1,c3-mag(eps4)+1)+1
        wpterm[k] = max(c1,mag(L+2)+c3-mag(eps3)+1)+1


# computing the tcoef[k,ell]
tcoef = {}
for k in range(0,L):
    for ell in range(0,3*k/2+1):
        tcoef[k,ell]=0

for k in range(0,L):
    mp.prec = wptcoef[k]
    for ell in range(0,3*k/2+1):
        tcoef[k,ell]=d[k,ell]*Fp[3*k-2*ell]/pipower[2*k-ell]
        tcoef[k,ell]=tcoef[k,ell]/((2*j)**ell)
\end{Verbatim}
}

\begin{proof}[Proof of the correctness of the algorithm]
It is clear that after lines 17--21 \(observe that the operation in
line 21 is free and for nothing\)  we get
\begin{displaymath}
\texttt{tcoef[k,ell]}
=\frac{(d^{(k)}_\ell+\alpha_{k,\ell})(F^{(3k-2\ell)}(p)
+\beta_{k,\ell})(1+\eta_3)}
{(2i)^\ell\pi^{2k-\ell}(1+\eta_{4k-2\ell})(1+\eta_1)}(1+\eta_1)
\end{displaymath}
where $|\alpha_{k,\ell}|\le \gamma_{k,\ell}$ and
$|\beta_{k,\ell}|\le\widetilde{\varepsilon_5}(m)$, with
$m=3k-2\ell$.

We may apply Lemmas \ref{lemma2} and \ref{lemma3} so that
\begin{displaymath}
\texttt{tcoef[k,ell]}
=(2i)^{-\ell}\pi^{\ell-2k}(d^{(k)}_\ell+\alpha_{k,\ell})
(F^{(3k-2\ell)}(p)+\beta_{k,\ell})(1+\eta_{4k-2\ell+6}).
\end{displaymath}
\noindent Hence, the absolute value of the difference $\texttt{tcoef[k,ell]}
-t^{(k)}_\ell $ is bounded by
\begin{multline*}
|\eta_{4k-2\ell+6}|2^{-\ell}\pi^{\ell-2k}(D^{(k)}_\ell+\gamma_{k\ell})(F_{m}
+\widetilde{\varepsilon_5}(m))+\\+
2^{-\ell}\pi^{\ell-2k}(D^{(k)}_\ell\widetilde{\varepsilon_5}(m)+
\gamma_{k\ell}F_{m}+\gamma_{k\ell}\widetilde{\varepsilon_5}(m)).
\end{multline*}

By lines 10 and 5 we know that $\texttt{wptcoef[k]}\ge \texttt{c1}$
so that all our working precisions $d_k$ satisfy $40(L+2)
2^{-d_k}<1$. Since $4k-2\ell+6\le 4(L+2)$,  Lemma \ref{boundExp}
applies and $|\eta_{4k-2\ell+6}|\le 4.24\, (L+2)2^{-d_k}$, so that
$2|\eta_{4k-2\ell+6}|\le 2\times0.106<1$. So the error is bounded by
\[
4.24(L+2) 2^{-d_k}\cdot 2^{-\ell}\pi^{\ell-2k} D^{(k)}_\ell
F_{m}+ 2\cdot 2^{-\ell}\pi^{\ell-2k}(D^{(k)}_\ell
\varepsilon_5\widetilde{\varepsilon_5}(m)+
\gamma_{k\ell}F_{m}+\gamma_{k,\ell} \widetilde{\varepsilon_5}(m)).
\]
By \eqref{fordk} we have
\begin{displaymath}
4.24\;(L+2) 2^{-d_k}\cdot 2^{-\ell}\pi^{\ell-2k} D^{(k)}_\ell
F_{m}\le 4.24\;(L+2) 2^{-d_k}A(\sigma)
\frac{\Gamma(k+1/2)^{1/2}}{(B_1 \sqrt{\pi})^k}.
\end{displaymath}
This will be $< \frac{\delta_k}{4}=\frac{\varepsilon_4 a^k}{16}$ if
we take
\begin{equation}
68 \frac{(L +2) A(\sigma)}{\varepsilon_4}\frac{\Gamma(k+1/2)^{1/2}}{(B_1
a \sqrt{\pi})^k}\le 2^{d_k}.
\end{equation}
This inequality is true for the choice  $d_k=\texttt{wptcoef[k]}$
in lines 4--9.

Now the choice of $\widetilde{\varepsilon_5}(m)$ is made in such a
way that
\begin{displaymath}
2\cdot
2^{-\ell}\pi^{\ell-2k}D^{(k)}_\ell\widetilde{\varepsilon_5}(m) <
\frac{\delta_k}{4}=\frac{\varepsilon_4
a^k}{16}\quad\Longleftrightarrow\quad \widetilde{\varepsilon_5}(m)<
\frac{1}{32}\frac{a^k\varepsilon_4}{2^{-\ell}\pi^{\ell-2k}
D^{(k)}_j}.
\end{displaymath}
Since $32\pi^2=315.827<316$ we have by \eqref{ukj}
\begin{displaymath}
\frac{1}{32}\frac{a^k\varepsilon_4}{2^{-\ell}\pi^{\ell-2k}
D^{(k)}_j}= \frac{u_{k,\ell}}{32}\ge \frac{\varepsilon_4(\pi^2A
B_1)^{m/3}}{32\pi^2 A(\sigma)}
\Bigl(\frac{m!}{\Gamma(m/3+2)}\Bigr)^{1/2}>\widetilde{\varepsilon_5}(m).
\end{displaymath}

Choosing
\begin{equation}\label{defgammak}
\gamma_{k,\ell}=\frac{1}{4}\cdot\frac{\delta_k}{8\cdot2^{-\ell}\pi^{\ell-2k}F_m}=
\frac{\sqrt{2\pi}}{128}\Bigl(\frac{\pi
a^2}{8}\Bigr)^{k/2}\frac{2^{2\ell}\varepsilon_4}{\Gamma\bigl(\frac{3k-2\ell+1}{2}\bigr)}
\end{equation}
we will have
\begin{displaymath}
2\cdot
2^{-\ell}\pi^{\ell-2k}\gamma_{k\ell}F_{m}<\frac{\delta_k}{4}=\frac{\varepsilon_4
a^k}{16}.
\end{displaymath}

Then, by the definition of $\widetilde{\varepsilon_5}(m)$ in
\eqref{defeps5}, we have
\begin{displaymath}
2\cdot 2^{-\ell}\pi^{\ell-2k}\gamma_{k\ell}\varepsilon_5(m)\le
8\cdot 2^{-\ell}\pi^{\ell-2k}\gamma_{k\ell}F_m=\frac{\delta_4}{4}.\qedhere
\end{displaymath}
\end{proof}

\subsection{Technical Lemmas.}
In the sequel we will need some concrete  inequalities. We will collect them here.

\begin{lemma}
For all $x>0$ we have the following inequality
\begin{equation}\label{gammainequality}
1>
\frac{\Gamma(x+1)\Gamma(x/3+2)}{(x+1)^2\bigl(\frac{(x+1)^4}{3e^4}\bigr)^{x/3}}>
\frac{2\pi}{3\sqrt{3}e^{4/3}}\approx 0.318741
\end{equation}
\end{lemma}

\begin{lemma}
For $x\ge3$
\begin{equation}\label{secondinegamma}
1.05599\approx\frac{3^{7/6}\sqrt{\pi\Gamma(10/3)}}{2^{10/3}}
\ge \Bigl(\frac{3\Gamma(x/3+7/3)}{2\Gamma(x/3+2)}\Bigr)^{1/2}
\frac{\Gamma(x/2+1)}{\Gamma(x/2+3/2)}
\Bigl(\frac{x+1}{3}\Bigr)^{1/3}>1.
\end{equation}
\end{lemma}

\begin{lemma}\label{fincreasing}
The function
\begin{equation}
\frac{5^x }{(2\pi)^{\frac{x-1}{2}}\Gamma(x/2+1/2)}
\Bigl(\frac{\Gamma(x+1)}{\Gamma(x/3+2)}\Bigr)^{1/2}
\end{equation}
is increasing for $x>0$.
\end{lemma}

\begin{proof}
We leave the proofs to the reader.
\end{proof}

\subsection{Computing $F^{(m)}(p)$.  Reduction to a polynomial.}
Our problem now is to compute
$\texttt{Fp[m]}:=F^{(m)}(p)+\Turing(\widetilde{\varepsilon_5}(m))$
for $0\le m\le 3L-3$.

If
\begin{equation}\label{Fmcondition}
F_m<\varepsilon_5(m):=\frac{(\pi^2 B_1 a)^{m/3}
}{316\,A(\sigma)}\Bigl(\frac{m!}
{\Gamma(m/3+2)}\Bigr)^{1/2}\,\varepsilon_4
\end{equation}
then $F_m\le \widetilde{\varepsilon_5}(m)$, and we can take
$\texttt{Fp[m]}=0$. By Lemma \ref{fincreasing} we have that
$\varepsilon_5(m)/F_m$ is increasing for $\pi^2 B_1 a>5^3$, so for
$t>100$ that we are assuming.  Therefore the condition
\eqref{Fmcondition} is frequently satisfied for $M\le m\le 3L-3$,
for some $M$. So we will solve our problem putting
\begin{equation}
\texttt{Fp[m]}=0,\qquad M\le m\le 3L-3,
\end{equation}
and solving the problem
\begin{equation}
\texttt{Fp[m]}=F^{(m)}(p)+\Turing(\varepsilon_5(m)),\qquad 0\le m<
M.
\end{equation}
In fact for $m< M$, we have $4F_m>F_m\ge \varepsilon_5(m)$, so that
\begin{displaymath}
\widetilde{\varepsilon_5}(m)=\max\{4F_m,\varepsilon_5(m)\}=\varepsilon_5(m).
\end{displaymath}

We
will compute the derivatives by means of Taylor series. Our first objective
is to reduce the Taylor series to a polynomial.  We recall that the
Taylor series for $F$ is given by (see  \cite{A86}*{formulas (47) and (56)})
\begin{displaymath}
F(z)=\sum_{j=0}^\infty c_{2j} z^{2j}, \qquad |c_{2j}|\le \frac{\pi^j}{2^{j+1}\,j!}.
\end{displaymath}

We shall use a Taylor polynomial $P(z)=\sum_{j=0}^{J-1} c_{2j}
z^{2j}$ to approximate $F(z)$.

First we define $J$ as the first   natural  number $J\ge 12$  such
that
\begin{equation}\label{defJ}
\frac{(2\pi)^J}{J!}\le \frac{\varepsilon_4}{632\,A(\sigma)},\qquad
\frac{(2\pi)^J}{J!}\le
\frac{\varepsilon_4}{632\,A(\sigma)}\frac{1}{M}\Bigl(\frac{\pi^2 B_1
a \sqrt{3} e^2} {M^2}\Bigr)^{(M-1)/3}.
\end{equation}

\begin{proposition}\label{PropAprTaylor}
Let $J$ be the natural number defined in \eqref{defJ} then
\begin{equation}\label{midapprox}
\Bigl|F^{(m)}(p)- P^{(m)}(p)\Bigr|\le
\frac{\varepsilon_5(m)}{2},\qquad (0\le m< M)
\end{equation}
where $P(z)$ is the polynomial $P(z)=\sum_{j=0}^{J-1} c_{2j} z^{2j}$
\end{proposition}

\begin{proof}
The rest of the series for $\frac{1}{m!}F^{(m)}(p)$ is bounded as follows \label{trece}
\begin{multline*}
\frac{1}{m!}\Bigl|\sum_{j\ge J}(2j)(2j-1)\cdots (2j-m+1) c_{2j} p^{2j-m}\Bigr|\le
\sum_{j\ge J}2^{2j}|c_{2j}|\le\\
\le \frac{1}{2}\sum_{j\ge J}\frac{(2\pi)^j}{j!}=
\frac{1}{2}\frac{(2\pi)^J}{J!}\Bigl(1+\frac{2\pi}{J+1}+
\frac{(2\pi)^2}{(J+1)(J+2)}+\cdots\Bigr)<
\end{multline*}
\(assuming $J+1\ge 13>4\pi$\)
\begin{displaymath}
<\frac{1}{2}\frac{(2\pi)^J}{J!}\sum_{j=0}^\infty \frac{1}{2^j}=\frac{(2\pi)^J}{J!}.
\end{displaymath}
Hence we must choose $P$ \(that is $J$\) in such a way that
\begin{equation}\label{aproxFP}
|F^{(m)}(p)-P^{(m)}(p)|< m!\frac{(2\pi)^J}{J!}\le
\frac{\varepsilon_5(m)}{2}.
\end{equation}

Since we want the polynomial $P$ \(and so $J$\) to be independent of
$m$ we want that the above inequality to be true for all $0\le m<
M$. By  definition \eqref{Fmcondition} we need to choose $J$ in
such a way that
\begin{equation}\label{secondcondition}
\frac{(2\pi)^J}{J!} \le \frac{(\pi^2 B_1 a)^{m/3}\varepsilon_4}{632
A(\sigma)} \frac{1}{\sqrt{m!\,\Gamma(m/3+2)}},\qquad 0\le m< M.
\end{equation}

Observe that the logarithm of the right hand side in
\eqref{secondcondition} as a function of $m$ is logarithmically
concave, so that the inequality above will be true for all $0\le m<
M$ if and only if it is true for the extremes $m=0$ and $m= M-1$. By
the inequality \ref{gammainequality} we have
\begin{multline*}
\frac{\varepsilon_5(m)}{2\, m!}=\frac{(\pi^2 B_1
a)^{m/3}\varepsilon_4}{632 A(\sigma)}
\frac{1}{\sqrt{m!\,\Gamma(m/3+2)}}>\\
> \frac{(\pi^2 B_1 a)^{m/3}\varepsilon_4}{632
A(\sigma)}\frac{1}{(m+1)\bigl(\frac{(m+1)^2}{\sqrt{3}e^2}\bigr)^{m/3}}=\\
=\frac{\varepsilon_4}{632
A(\sigma)}\frac{1}{(m+1)}\Bigl(\frac{\pi^2 B_1 a
\sqrt{3} e^2} {(m+1)^2}\Bigr)^{m/3}.
\end{multline*}
Hence, the last two conditions on $J$ given in \eqref{defJ} imply
\eqref{aproxFP} for $0\le m<M$.
\end{proof}

\subsection{Computing the derivatives of a polynomial.}\label{SectionPolynomial}

In this Section we consider the general problem of computing the derivatives of a polynomial
that we will apply in Section \ref{sectionFmp} to our particular case.

We assume that $P(x)=\sum_{k=0}^{J-1} c_k x^k$ is a given
polynomial. We want to compute the derivatives
$P^{(m)}(p)+\Turing(\varepsilon_m)$ for $0\le m<M$ at a real point
$p$ with $|p|\le1$. We assume that the coefficients $c_k$ are
difficult to compute, so we want  to compute
$\texttt{c[k]}=c_k+\Turing(\delta_k)$ only once. Also we want to operate in a
fixed working precision when computing the polynomials. Our problem
is to determine what must be the values of the $\delta_k$ and what
will be the adequate working precision.

\begin{proposition}\label{PropPolynomial}
In order to  compute the derivatives
$\texttt{Pp[m]}:=P^{(m)}(p)+\Turing(\varepsilon_m)$ for $0\le m<M$
we choose $\delta_k\ge0$ such that
\begin{equation}\label{deltacondition}
\sum_{k=m}^{J-1}2^k
\delta_k<\frac{1}{3}\frac{\varepsilon_m}{m!},\qquad 0\le m<M.
\end{equation}
Then compute the representable numbers
$\texttt{c[k]}=c_k+\Turing(\delta_k)$.\footnote{$\delta_k=0$ may be
used only in case $c_k$ is representable. In this case
$\texttt{c[k]}=c_k$.}

Choose the working precision $d$ as the least natural number such that
\begin{equation}\label{wpcomputPol}
2^d>22 J,\qquad  2^d>6.36 J \frac{m! H_m}{\varepsilon_m},\qquad 0\le
m<M
\end{equation}
where
\begin{equation}\label{defcoefAm}
H_m:=\sum_{k=m}^{J-1} 2^k |c_k|.
\end{equation}
Then we get \texttt{Pp[m]} for $0\le m< M$ by the following
procedure.
\end{proposition}

\bigskip
\bigskip
{\footnotesize
\begin{Verbatim}[numbers= left,numbersep=2pt,numberblanklines=false,frame=single,
label=\fbox{\Large Computing \texttt{Pp[m]}.}]

mp.prec = d

for m in range(0,M):
    sumP = 0
    for k in range(J-m-1,-1,-1):
        sumP = (sumP * p)+ c[k]
    Pp[m] = sumP
    # preparation of the new coefficients
    for k in range(0,J-m-1):
        c[k] = (k+1)* c[k+1]
\end{Verbatim}
}

\begin{proof}
Let $\texttt{c[k]}=c_k+\alpha_k$ with $|\alpha_k|\le \delta_k$, and consider the polynomial \[Q(x)=\sum_{k=0}^{J-1} \texttt{c[k]} x^k.\] Then we have
\begin{multline*}
|P^{(m)}(p)-Q^{(m)}(p)|=\Bigl|\sum_{k=m}^{J-1}k(k-1)\cdots(k-m+1)\alpha_k p^{k-m}\Bigr|\le\\
\le \sum_{k=m}^{J-1}k(k-1)\cdots(k-m+1) \delta_k.
\end{multline*}
It follows that
\begin{equation} \label{menor13}
\frac{1}{m!}|P^{(m)}(p)-Q^{(m)}(p)|\le \sum_{k=m}^{J-1} 2^k \delta_k<\frac13\frac{\varepsilon_m}{m!}
\end{equation}
by the choice of the $\delta_k$

The coefficients $\texttt{c[k]}$ change in each run of the for loop in $m$ \(see lines 8--9\).  In the first run \($m=0$\) we compute with a rounded version of these  numbers $\texttt{c[k]}(1+\eta_1)$. In the consecutive runs  we change the numbers to the \(rounded\)
\begin{align*}
m=1\qquad &\texttt{c[k]} = (k+1)\texttt{c[k+1]}(1+\eta_2)\\
m=2\qquad &\texttt{c[k]} = (k+1)(k+2)\texttt{c[k+2]}(1+\eta_3)\\
\dots\\
m\qquad & \texttt{c[k]} = (k+1)(k+2)\cdots (k+m)\texttt{c[k+m]}(1+\eta_{m+1})
\end{align*}
where we have put the numbers $\texttt{c[k]}$ in the $m$--th run in terms of the initial values of the $\texttt{c[k]}$.

In the $m$--th run of the first for loop,   the consecutive values taken by the variable \texttt{sumP} \(see line 5\)  are:
\begin{align*}
k=J-m-1\qquad & \texttt{sumP} = \texttt{c[J-m-1]}\\
k=J-m-2\qquad & \texttt{sumP} = \texttt{c[J-m-1]p}(1+\eta_2)+c[J-m-2](1+\eta_1)\\
k=J-m-3\qquad & \texttt{sumP} = \texttt{c[J-m-1]p}^2(1+\eta_4)+\texttt{c[J-m-2]p}(1+\eta_3)+\\
&\hspace{5.5cm}+\texttt{c[J-m-3]}(1+\eta_1)\\
k=0\qquad & \texttt{sumP} = \texttt{c[J-m-1]p}^{J-m-1}(1+\eta_{2J-2m-2})+\\
&\hspace{4cm}+\sum_{k=0}^{J-m-2}\texttt{c[k]p}^k(1+\eta_{2k+1}).
\end{align*}
Hence, in line   6  we get \(changing the first factor
$(1+\eta_{2J-2m-2})$ to $(1+\eta_{2J-2m-1})$\)
\begin{displaymath}
\texttt{Pp[m]}=\sum_{k=0}^{J-m-1}\texttt{c[k]p}^k (1+\eta_{2k+1}).
\end{displaymath}
Then in terms of the initial \texttt{c[k]}
\begin{displaymath}
\texttt{Pp[m]}=\sum_{k=0}^{J-m-1}(k+1)(k+2)\cdots(k+m)\texttt{c[k+m]p}^k
(1+\eta_{2m+2k+2}).
\end{displaymath}
\(We put here $2m$ instead of $m$ to simplify the reasoning below.\)

It follows that
\begin{displaymath}
\texttt{Pp[m]}=Q^{(m)}(p)+R_m
\end{displaymath}
where
\begin{displaymath}
|R_m|\le
\sum_{k=0}^{J-m-1}(k+1)(k+2)\cdots(k+m)|\texttt{c[k+m]}|\cdot
|\eta_{2m+2k+2}|.
\end{displaymath}
Hence
\begin{multline*}
\frac{1}{m!}|Q^{(m)}(p)-\texttt{Pp[m]}|\le
\sum_{k=0}^{J-m-1}2^{k+m}|\texttt{c[k+m]}|\cdot |\eta_{2m+2k+2}|\le \\
\le \sum_{k=0}^{J-m-1}2^{k+m}(|c_{k+m}|+\delta_{k+m})\,
|\eta_{2m+2k+2}|=\sum_{k=m}^{J-1} 2^k (|c_k|+\delta_k)|\eta_{2k+2}|.
\end{multline*}
Our choice  of the working precision $d$ satisfies $2J\, 2^{-d}<0.1$
so that Lemma \ref{boundExp} applies and
\begin{displaymath}
\frac{1}{m!}|Q^{(m)}(p)-\texttt{Pp[m]}|\le 2.12\,J\, 2^{-d}
\sum_{k=m}^{J-1} 2^k (|c_k|+\delta_k)= 2.12\,J\, 2^{-d}
\Bigl(H_m+\frac13\frac{\varepsilon_m}{m!}\Bigr)
\end{displaymath}
with $H_m$ defined in \eqref{defcoefAm}. Our choice of $d$
in \eqref{wpcomputPol} guarantees that
\begin{equation}\label{menor23}
\frac{1}{m!}|Q^{(m)}(p)-\texttt{Pp[m]}|< \frac23\frac{\varepsilon_m}{m!}.
\end{equation}
Finally \eqref{menor13} and \eqref{menor23} prove our Theorem.
\end{proof}

\subsection{Computing \texttt{Fp[m]}.} \label{sectionFmp}
After \eqref{midapprox} our problem is to compute $P^{(m)}(p)+\Turing(\varepsilon_5(m)/2)$, where $P$ is the Taylor polynomial considered in Proposition \ref{PropAprTaylor} and  $J$ is defined in  \eqref{defJ}. We may apply Proposition \ref{PropPolynomial} to solve our problem.

\begin{proposition}\label{PropComputPolynomial}
Assume that  condition \eqref{thirdcondition} is satisfied.
To get $P^{(m)}(p)+\Turing(\varepsilon_5(m)/2)$ we first
compute
\begin{equation}\label{choosedeltak}
c_{2k}+\Turing(2^{-2k}\varepsilon_6),\qquad (0\le k\le J)\quad
\text{where}\quad
\varepsilon_6:=\frac{1}{3J}\frac{(2\pi)^J}{J!}.
\end{equation}
We fix the working precision \texttt{wpfp} as the least natural number $d$ such that
\begin{equation}\label{wpfmp}
2^d>44 \,J,\qquad 2 ^d > 6812\, J \frac{m!}{\varepsilon_5(m)},\qquad 0\le m\le \min(3L-3,21)
\end{equation}
and then follow the procedure explained in Section \ref{SectionPolynomial}.
\end{proposition}

\begin{proof}
Our polynomial $P(x)=\sum_{k=0}^{J-1} c_{2k}x^{2k}$ has the odd numbered coefficients equal to $0$, so that we can take $\texttt{c[2k+1]}=0$ and therefore $\delta_{2k+1}=0$. Also we choose $\delta_{2k}= 2^{-2k}\varepsilon_6$ as indicated in \eqref{choosedeltak}.  With these choices the sums in \eqref{deltacondition} are
\begin{displaymath}
\sum_{k=m}^{2J-2}2^k \delta_k<\sum_{m/2 \le k\le J-1}2^{2k}2^{-2k}\varepsilon_6\le J \varepsilon_6\le \frac{1}{3}\frac{(2\pi)^J}{J!}.
\end{displaymath}
Since we have chosen $J$ \(see \eqref{aproxFP}\) in such a way that $\frac{(2\pi)^J}{J!}\le
\frac{\varepsilon_5(m)}{2\, m!}$, the condition \eqref{deltacondition} of Proposition \ref{PropPolynomial} is satisfied.

Now we must determine in our case the working precision of Proposition \ref{PropPolynomial}  and defined in \eqref{wpcomputPol}.

In our case the numbers $H_m$ are given by
\begin{displaymath}
H_{2m-1}=H_{2m}=\sum_{k=m}^{J-1} |c_{2k}| 2^{2k}\le \frac{1}{2}\sum_{k=m}^{\infty}\frac{(2\pi)^k}{k!}
\end{displaymath}
\(after the bound  of $|c_{2k}|$ given in \cite{A86}*{equation (56)}\).

In the proof of Proposition \ref{PropAprTaylor}  we have seen that
$H_{2m}\le(2\pi)^m/m!$ when $m\ge12$. One may verify that this also holds for $m\ge11$.

Therefore, for even $k\ge 22$ we have $H_k\le
(2\pi)^{k/2}/\Gamma(k/2+1)$. For $k=2m-1$  odd we have
$H_{k}=H_{2m}\le
(2\pi)^{m}/\Gamma(m+1)<(2\pi)^{m-1/2}/\Gamma(m+1/2)$
 for $m\ge11$, so that
\begin{displaymath}
H_m\le \frac{(2\pi)^{m/2}}{\Gamma(m/2+1)},\qquad (m\ge21).
\end{displaymath}

By \eqref{wpcomputPol} we may take as working precision the
 least natural number $d$ such that
\begin{equation}\label{dcondfirstform}
2^d>44 \,J,\qquad 2 ^d > 12.72\, J \frac{m!
H_m}{\varepsilon_5(m)/2}, \qquad 0\le m< M.
\end{equation}

For $m\ge21$ the last condition follows from
\begin{displaymath}
2 ^d > 25.44\, J \frac{m! }{\varepsilon_5(m)}\frac{(2\pi)^{m/2}}{\Gamma(m/2+1)}:=K(m),\qquad 21\le m\le 3L-3.
\end{displaymath}

We have  $K(m+1)<K(m)$. In fact applying \eqref{secondinegamma}
\begin{multline*}
\frac{K(m+1)}{K(m)}=\frac{\sqrt{2\pi(m+1)}}{(\pi^2B_1a)^{1/3}}\Bigl(
\frac{\Gamma(m/3+7/3)}{\Gamma(m/3+2)}\Bigr)^{1/2}\frac{\Gamma(m/2+1)}{\Gamma(m/2+3/2)}\le\\
\le 1.056\frac{\sqrt{2\pi(m+1)}}{(\pi^2B_1a)^{1/3}}\cdot \frac{\sqrt{2}}{\sqrt{3}}\Bigl(\frac{3}{m+1}\Bigr)^{1/3}<
3.12\, \Bigl(\frac{\sqrt{m+1}}{\pi^2 B_1 a}\Bigr)^{1/3}
\le 3.12\, \Bigl(\frac{\sqrt{M}}{\pi^2 B_1 a}\Bigr)^{1/3}.
\end{multline*}
Since $M\le 3L-3$ and we assume that  $3L<2a^2/25$ \(see
\eqref{thirdcondition}\) and by definition $B_1\ge1$ we get
\begin{displaymath}
\frac{K(m+1)}{K(m)}<
3.12\,\Bigl(\frac{\sqrt{2}}{5\pi^2}\Bigr)^{1/3}\le0.95478<1.
\end{displaymath}

Therefore, the conditions on $d$ given in \eqref{dcondfirstform} are
equivalent for $M\ge 21$ to
\begin{displaymath}
2^d>44 \,J,\quad 2 ^d > 25.44\, J \frac{m! \,
H_m}{\varepsilon_5(m)},\quad 0\le m\le 20,\quad 2^d> 12.72\, J
\frac{21!\, H_{21}}{\varepsilon_5(21)/2}
\end{displaymath}
and in case  $M\le 20$ it will be
\begin{displaymath}
2^d>44 \,J,\qquad 2 ^d > 25.44\, J \frac{m!
H_m}{\varepsilon_5(m)},\qquad 0\le m\le M.
\end{displaymath}

It is also convenient to observe that $H_m\le 267.746$, so that we can put our conditions in the following form:
\begin{equation}
2^d>44 \,J,\qquad 2 ^d > 6812\, J \frac{m!}{\varepsilon_5(m)},\qquad
0\le m\le \min(M,21).\qedhere
\end{equation}
\end{proof}

\subsection{Computing the coefficients $\texttt{c[2n]}=c_{2n}+\Turing(\delta_{2k})$.}
The $c_{2k}$ are the coefficients of the Taylor expansion of the entire function $F(z)=\sum_{n=0}^\infty c_{2n} z^{2n}$. In \cite{A86}*{(48)} we get the formula
\begin{multline}\label{c2n}
c_{2n}=-\frac{i}{\sqrt{2}}\Bigl(\frac{\pi}{2}\Bigr)^{2n}\sum_{k=0}^n\frac{(-1)^k}{(2k)!}
2^{2n-2k}\frac{(-1)^{n-k}E_{2n-2k}}{(2n-2k)!}+\\+e^{3\pi
i/8}\sum_{j=0}^n(-1)^j\frac{
E_{2j}}{(2j)!}\frac{i^{n-j}\pi^{n+j}}{(n-j)!2^{n-j+1}}
\end{multline}
which gives $c_{2n}$ in terms of the Euler numbers.

The number $c_{2n}$ is a polynomial of degree $2n$ in $\pi$  and coefficients
in $\Q(e^{3\pi i/8})$, so this expression can not be simplified.

Since $F$ is an entire function, for all $R>0$, we have $c_{2n}
R^{2n}\to0$ when $n\to\infty$, but each of the two summands in
\eqref{c2n} is of the order $\frac{2}{\pi} 2^{2n}$.   Hence the two
sums in \eqref{c2n} cancel out very much,
  so that computing
these coefficients to a given precision is a time consuming task.
For this reason these data once computed will be stored in the
cache.

It is easy to see that
\begin{equation}
c_{2n}=(-1)^{n+1}\frac{i}{\sqrt{2}}\sum_{k=0}^n (-1)^k V(k)W(2n-2k)+\frac{e^{3\pi i /8}}{2}
\sum_{k=0}^n i^{n-k}V(k)W(n-k)
\end{equation}
where we have put
\begin{equation}\label{defVW}
V_n=(-1)^n E_{2n}\frac{\pi^{2n}}{(2n)!},\qquad
W_n=\frac{1}{n!}\Bigl(\frac{\pi}{2}\Bigr)^n.
\end{equation}
Hence we define the two scalar products
\begin{align}
P_1(n):&=(-1)^{n+1} i \sum_{k=0}^n (-1)^k V_kW_{2n-2k}\label{PP1}\\
P_2(n):&=\sum_{k=0}^n i^{n-k}V_kW_{n-k}.\label{PP2}
\end{align}
Then $c_{2n}= \mu P_1(n)+\nu P_2(n)$, where $\mu=2^{-1/2}$ and $\nu=e^{3\pi i/8}/2$.

The Euler numbers satisfy
\begin{equation}
(-1)^n\frac{E_{2n}}{(2n)!}\left(\frac{\pi}{2}\right)^{2n}=\frac{4}{\pi}\sum_{k=0}^\infty
\frac{(-1)^k}{(2k+1)^{2n+1}}
\end{equation}
so that
\begin{equation}\label{boundVW}
|V_n|\le \frac{4}{\pi}2^{2n},\qquad |W_n|\le \frac{1}{n!}\Bigl(\frac{\pi}{2}\Bigr)^n.
\end{equation}
From \eqref{PP1}, \eqref{PP2} and \eqref{boundVW} it follows that
\begin{equation}\label{boundP1P2}
|P_1(n)|< \frac{4}{\pi}\cosh(\pi/4) 2^{2n}<1.69\cdot2^{2n}, \quad
|P_2(n)|<\frac{4}{\pi}e^{\pi/8} 2^{2n}<1.89\cdot2^{2n}.
\end{equation}

\begin{proposition}
Assume that we have solved the two problems:
\begin{equation}\label{problemP1P2}
\texttt{P1[n]}= P_1(n)+\Turing(2^{-2n}\varepsilon_6/6),\qquad
\texttt{P2[n]}= P_2(n)+\Turing(2^{-2n}\varepsilon_6/6).
\end{equation}
Then we get $\texttt{c[2n]}=c_{2n}+\Turing(\delta_{2n})$ by the following procedure:
\end{proposition}

\bigskip
\bigskip
{\footnotesize
\begin{Verbatim}[numbers= left,numbersep=2pt,numberblanklines=false,frame=single,
label=\fbox{\Large Computing \texttt{c[2n]}.}]

mp.prec=15
wpc0 = 5 - mag(eps6)
wpc = max(6,4*J+wpc0)
mp.prec = wpc
mu = sqrt(2)/2
nu = exp(3*pi*j/8)/2


c={}

for n in range(0,J):
    mp.prec = 15
    wpc = max(6,4*n+wpc0)
    mp.prec = wpc
    c[2n] = mu*P1[n] + nu*P2[n]
\end{Verbatim}
}

\begin{proof}
In lines 5--6 we compute the constants $\mu$ and $\nu$ with a precision  $d'$. We have  $d'\ge d+4$ for each value $d$ of the working precision used in line 12.

Let $d'$ be the value of the variable \texttt{wpc} in lines 4--6. Then we have
\begin{displaymath}
\texttt{mu} = 2^{-1/2}\,(1+\eta_1),\qquad \texttt{nu} =
\frac{e^{3\pi i/8(1+\eta_1)}}{2}(1+\eta_1),
\end{displaymath}
\(recall that division by 2 or multiplication by $j$ is for free and nothing, and our hypothesis that the $\exp$ and $\sqrtt$ are well implemented\).

Since $d'\ge2$ we have $e^{3\pi i/8(1+\eta_1)}=e^{3\pi i/8} e^{3\pi i\eta_1/8}$ and
\begin{displaymath}
|e^{3\pi i\eta_1/8} -1|\le \sum_{k=1}^\infty \Bigl(\frac{3\pi}{8}|\eta_1|\Bigr)^k\frac{1}{k!}\le \frac{3\pi}{8}|\eta_1|e^{3\pi /32}\le 1.582 \,|\eta_1|
\end{displaymath}
so that  $\texttt{nu} =\frac12 e^{3\pi i/8}(1+2\eta_1)(1+\eta_1)$.
When we change the working precision to one of the values $d$ we
will have $\texttt{nu} =\frac12 e^{3\pi i/8}(1+\eta_1)$, in fact it
is easy to see that if $\xi=(1+2\eta_1)(1+\eta_1)-1$ with
$|\eta|<2^{-d-4}$ then $|\xi|\le 2^{-d}$. Since we have to truncate
these two values, we in fact operate \(in line 12\) with rounded
numbers such that
\begin{displaymath}
\texttt{mu} = 2^{-1/2}\,(1+\eta_2)=\mu(1+\eta_2),\qquad \texttt{nu}
= \frac{e^{3\pi i/8}}{2}(1+\eta_2)=\nu(1+\eta_2).
\end{displaymath}

\texttt{mu} is real but \texttt{nu} is complex so that from line 12 we get
\begin{displaymath}
\texttt{c[2n]}=\mu(P_1(n)+\alpha)(1+\eta_5)+\nu(P_2(n)+\beta)(1+\eta_5)
\end{displaymath}
where $|\alpha|$, and $|\beta|\le \frac{1}{6}2^{-2n}\varepsilon_6$.

Hence $\texttt{c[2n]}=c_{2n}+R_n$ where
\begin{displaymath}
|R_n|\le
2^{-1/2}|P_1(n)\eta_5|+2^{-1/2}|\alpha(1+\eta_5)|+\frac{1}{2}
|P_2(n)\eta_5|+\frac{1}{2}|\beta(1+\eta_5)|.
\end{displaymath}
Since $d\ge6$ we have $5\cdot 2^{-d}<0.1$ and Lemma \ref{boundExp} applies, we also apply \eqref{boundP1P2} so that
\begin{multline*}
|R_n|\le 2^{-1/2}\cdot1.69\cdot2^{2n}\cdot5.3 \cdot 2^{-d}+\frac{1}{6\sqrt{2}}2^{-2n}\varepsilon_6\cdot1.09
+\frac{1}{2} 1.89\cdot 2^{2n}\cdot 5.3 \cdot2^{-d}+\\+
\frac{1}{12}2^{-2n}\varepsilon_6\cdot 1.09
<11.3421\cdot 2^{2n-d}+0.2193\cdot 2^{-2n}\varepsilon_6.
\end{multline*}
By the choice of $d$ we have $2^d>16\frac{2^{4n}}{\varepsilon_6}$ so
that
\begin{displaymath}
|R_n|<\frac{11.3421}{16}\cdot 2^{-2n}\varepsilon_6+0.2193\cdot 2^{-2n}\varepsilon_6<2^{-2n}\varepsilon_6=\delta_{2n}
\end{displaymath}
as was required in \eqref{choosedeltak}.
\end{proof}

\subsection{A Convolution.}\label{SectionConvolution}
We have to solve the problems \eqref{problemP1P2}. Recall that
$P_1(n)$ and $P_2(n)$ are given as convolutions \eqref{PP1}
and \eqref{PP2}. So in this Section we treat the general
problem and then we shall apply it to these two cases.

We have real numbers $V_k$ and $W_k$  $0\le k\le n$, with
$|V_k|\le A_k$ and $|W_k|\le B_k$.
We want to compute
$\sum_{k=0}^nV_kW_{n-k}+\Turing(\varepsilon)$.

Our problem is to determine a working precision \texttt{wp}
 and numbers $\mu_k$ and $\nu_k$ so that computing
\begin{equation}
\texttt{v[k]} = V_k+\Turing(\mu_k),\qquad \texttt{w[k]} = W_k+\Turing(\nu_k)
\end{equation}
and operating to the working precision \texttt{wp} we get the desired result.

Since we want to compute the convolutions for different values of $n$, and desire to compute the $V_k$ and $W_k$ only once, we require that $\mu_k$ and $\nu_k$ does not depend on $n$.

\begin{proposition}\label{PropConvolution}
In order to compute the convolution $\sum_{k=0}^nV_kW_{n-k}+\Turing(\varepsilon)$
we take the working precision \texttt{wp} as the least natural  number $d$ such that
\begin{equation}
2^d >10(n+4),\qquad 2^d >(2.385 M_n+2.5\,N_n)/\varepsilon
\end{equation}
where
\begin{equation}\label{defMN}
M_n=\sum_{k=0}^n (n-k+4)A_kB_{n-k},\quad
N_n=\sum_{k=0}^n A_kB_{n-k}
\end{equation}
and put $\mu_k= 2^{-d} A_k$ and  $\nu_k=2^{-d}B_k$,
then our problem can be solved by the following procedure:
\end{proposition}

\bigskip
\bigskip
{\footnotesize
\begin{Verbatim}[numbers= left,numbersep=2pt,numberblanklines=false,frame=single,
label=\fbox{\Large Computing $\texttt{P}=\sum_{k=0}^nV_kW_{n-k}+\Turing(\varepsilon)$.}]

mp.prec = wp

P = 0
for k in range(0,n+1):
    term= v[k]*w[n-k]
    P += term

\end{Verbatim}
}

\begin{proof}
The  values of the variable \texttt{term} are rounded numbers. Taking in consideration   the truncation and the error in the product
\begin{displaymath}
\texttt{term}=\texttt{v[k]}(1+\eta_1)\texttt{w[n-k]}(1+\eta_1)(1+\eta_1)=
(V_k+\alpha_k)(W_{n-k}+\beta_{n-k})(1+\eta_3)
\end{displaymath}
where $|\alpha_k|\le\mu_k$ and  $|\beta_k|\le\nu_k$.

The final value of the variable $\texttt{P}$ will be by Proposition \ref{sumofrounded} is
\begin{equation}
\texttt{P}=\sum_{k=0}^n (V_k+\alpha_k)(W_{n-k}+\beta_{n-k})(1+\eta_3)(1+\eta_{n-k+1})
\end{equation}
It follows that
\begin{multline*}
\Bigl|\texttt{P}-\sum_{k=0}^n V_kW_{n-k}\Bigr|\le\\ \sum_{k=0}^n (|V_k|+\mu_k)(|W_{n-k}|+\nu_{n-k})|\eta_{n-k+4}|+
\sum_{k=0}^n(|V_k|\nu_{n-k}+|W_{n-k}|\mu_k+\mu_k\nu_{n-k})
\end{multline*}
By the first condition imposed to the working precision
 $d$ we may apply in all cases Lemma \ref{boundExp}. The first sum contain  the terms
\begin{displaymath}
\sum_{k=0}^n |V_k W_{n-k}||\eta_{n-k+4}|\le 1.06\,2^{-d}\sum_{k=0}^n (n-k+4)A_kB_{n-k}=
1.06\,2^{-d} M_n
\end{displaymath}
where $M_n$ is defined in \eqref{defMN}.

With our choice of $\mu_k$ and $\nu_k$ we get
\begin{multline*}
\Bigl|\texttt{P}-\sum_{k=0}^n V_kW_{n-k}\Bigr|\le 1.06\,M_n 2^{-d}+
1.06\,2^{-d}\Bigl\{\sum_{k=0}^n(n-k+4)A_k 2^{-d}B_{n-k}+\\+\sum_{k=0}^n (n-k+4)B_k  2^{-d}A_{n-k}+
\sum_{k=0}^n (n-k+4)2^{-d}A_k2^{-d}B_{n-k}\Bigr\}+\\ +
\sum_{k=0}^n A_k2^{-d}B_{n-k}+\sum_{k=0}^nB_{n-k}2^{-d}A_{k}+\sum_{k=0}^n 2^{-d}A_k
2^{-d}B_{n-k}.
\end{multline*}
This can be written as
\begin{multline*}
\Bigl|\texttt{P}-\sum_{j=0}^n V_jW_{n-j}\Bigr|\le
1.06(1+2^{-d})^2M_n2^{-d}+(2+2^{-d})2^{-d}N_n\le\\ \le (2.385
M_n+2.5 N_n)2^{-d}<\varepsilon
\end{multline*}
by our choice of the working precision $d$.
\end{proof}

\subsection{Computing $P_1(n)+\Turing(2^{-2n}\varepsilon_6/6)$
and $P_2(n)+\Turing(2^{-2n}\varepsilon_6/6)$.}

Note that in both cases  we always have to multiply real numbers or
a real number and a purely imaginary. So we can use the most simple
bounds and the results of the analysis of Section
\ref{SectionConvolution} are applicable.

For $P_1(n)$ with the notations in Proposition \ref{PropConvolution} we have
\begin{align}
M_n&=\sum_{j=0}^n (n-j+4) \frac{4}{\pi} 2^{2j} \frac{1}{(2n-2j)!}
\left(\frac{\pi}{2}\right)^{2n-2j}
<\\
&\qquad <2^{2n}\Bigl(\frac{16}{\pi}\cosh(\pi/4)+\frac12
\sinh(\pi/4)\Bigr)<7.2 \cdot 2^{2n},\\
N_n&=\sum_{j=0}^n  \frac{4}{\pi} 2^{2j} \frac{1}{(2n-2j)!}
\left(\frac{\pi}{2}\right)^{2n-2j}
< 2^{2n}\frac{4}{\pi}\cosh(\pi/4)
<1.7\cdot 2^{2n}.
\end{align}
So following Proposition \ref{PropConvolution} we must operate with working precision \texttt{wpp1} defined as the least
natural number $d$ such that
\begin{equation}
2^d>10(n+4)\quad \text{and} \quad2^d\ge
128.6\cdot2^{4n}/\varepsilon_6 >6(2.385 M_n+2.5
N_n)/(\varepsilon_62^{-2n})
\end{equation}
and compute with the values obtained for
\begin{equation}
\texttt{v[k]}=V(k)+\Turing(2^{-d} \frac{4}{\pi} 2^{2k} ),\quad
\texttt{w[2k]}=W(2k) +\Turing\Bigl(2^{-d}\frac{1}{(2k)!}\left(\frac{\pi}{2}\right)^{2k}\Bigr),
\end{equation}
for $0\le k\le n$.

For $P_2(n)$ we apply also Proposition \ref{PropConvolution}. In this case
\begin{align}
M_n&=\sum_{j=0}^n (n-j+4) \frac{4}{\pi} 2^{2j}
\frac{1}{(n-j)!}\left(\frac{\pi}{2}\right)^{n-j}
<\Bigl(\frac{16}{\pi}+\frac12\Bigr)e^{\pi/8} 2^{2n}<8.3
\cdot 2^{2n},\\
N_n&=\sum_{j=0}^n  \frac{4}{\pi} 2^{2j}
\frac{1}{(n-j)!}\left(\frac{\pi}{2}\right)^{n-j}
<\frac{4}{\pi} e^{\pi/8}2^{2n}
<1.9\cdot 2^{2n}.
\end{align}
Hence we must take the working precision \texttt{wpp2}
as the least natural number $d$ such that
\begin{equation}
2^d>10(n+4)\quad \text{and} \quad2^d\ge
147.3\cdot2^{4n}/\varepsilon_6
\end{equation}
and compute the convolution with
\begin{equation}
\texttt{v[k]}=V(k)+\Turing(2^{-d} \frac{4}{\pi} 2^{2k} ),\quad
\texttt{w[k]}=W(k) +\Turing\Bigl(2^{-d}\frac{1}{k!}\left(\frac{\pi}{2}\right)^k\Bigr),
\end{equation}
for $0\le k\le n$.

We need to compute $c_{2n}$ for $0\le n<J$. To compute $c_{2n}$ we need $P_1(n)$ and $P_2(n)$. Now to compute $P_1(n)$ and $P_2(n)$ we need the numbers $V(k)$, $W(k)$ and $W(2k)$ for $0\le k\le n$.  Since we want to compute $V(k)$ and $W(k)$ only one time, we shall compute them to the biggest accuracy needed.  That is we put the working precision \texttt{wpvw} as the least natural number $d$ such that
\begin{equation}
2^d>10(J+3)\quad \text{and} \quad2^d\ge2^{4J+4}/\varepsilon_6>
147.3\cdot2^{4J-4}/\varepsilon_6
\end{equation}
and with this working precision $\texttt{wpvw}=d$  compute
\begin{equation}
\texttt{v[k]}=V(k)+\Turing(2^{-d} \frac{4}{\pi} 2^{2k} ),\quad
\texttt{w[k]}=W(k) +\Turing\Bigl(2^{-d}\frac{1}{k!}\left(\frac{\pi}{2}\right)^k\Bigr),
\end{equation}
for $0\le k<J$.

\subsection{Computing $d^{(n)}_k+\Turing(\gamma_{n,k})$.}

Recall the definition of these numbers
\begin{align}
&d^{(0)}_0=1, \quad d^{(0)}_k =0 \quad \text{for } k\ne0\\
&d^{(n)}_k=0 \quad \text{for }\quad k<0 \quad\text{and for }\quad k>3n/2\\
&2(3n-2k) d^{(n)}_k=\tfrac12 
d^{(n-1)}_k+(1-2\sigma)d^{(n-1)}_{k-1}-2(3n-2k)(3n-2k+1)
d^{(n-1)}_{k-2}\label{recurrenced}\\
&d^{(n)}_{3n/2} =-\sum_{k=0}^{3n/2 -1} (-1)^{3n/2-k} d^{(n)}_k
\frac{(3n-2k)!}{(3n/2-k)!},\qquad  3n \equiv 0\pmod
2.\label{case3n2k}
\end{align}

The first values of these coefficients are:
\begin{align}
d^{(1)}_0&=\frac{1}{12},\quad d^{(1)}_1=-(\sigma-\tfrac12 ),\\
d^{(2)}_0&=\frac{1}{288},\quad
d^{(2)}_1=-\frac{1}{12}(\sigma-\tfrac12 ),\quad d^{(2)}_2=\frac12(\sigma-\tfrac12 )^2-\frac14,\notag\\
&\hspace{3cm}d^{(2)}_3=(\sigma-\tfrac12 )^2+(\sigma-\tfrac12 )-\frac{1}{12}.
\end{align}
In general it can be shown that $d^{(n)}_k$ is a polynomial in
$\sigma$ of degree $\max(k,n)$.

\begin{proposition}\label{propd}
Assume that $|\sigma|\le \frac{a}{2}$. In order to compute
$\texttt{d[n,k]}=d^{(n)}_k+\Turing(\gamma_{n,k})$ for $0\le k\le
3n/2$ we define for each value of $n$ the working precision
$\texttt{wpd[n]}$ as the least natural number $d\ge6$ such that
\begin{equation}
2^d>40L^2,\qquad 2^d > (1+|\sigma|)\frac{A\,2^{13}}{
\varepsilon_4}\Bigl(\frac{8}{\pi a^2 B_1^2}\Bigr)^{n/2}
\Gamma(n-\tfrac12 )^{1/2}.
\end{equation}
and apply the following procedure:
\end{proposition}

\bigskip
\bigskip
{\footnotesize
\begin{Verbatim}[numbers= left,numbersep=2pt,numberblanklines=false,frame=single,
label=\fbox{\Large Computing $d^{(n)}_k$.}]

mp.prec = wpd[1]
psigma = 1-(2*sigma)

d = {}

d[0,-2]=0; d[0,-1]=0; d[0,0]=1; d[0,1]=0

for n in range(1,L):
    mp.prec = wpd[n]
    for k in range(0,3*n/2+1):
        m = 3*n-2*k
        if(m!=0):
            m1 = mpf('1')/m
            c1= m1/4
            c2=(psigma*m1)/2
            c3=-(m+1)
            d[n,k]=c3*d[n-1,k-2]+c1*d[n-1,k]+c2*d[n-1,k-1]
        else:
            d[n,k]=0
            for r in range(0,k):
                add=d[n,r]*(mpf('1.0')*fac(2*k-2*r)/fac(k-r))
                d[n,k] -=  ((-1)**(k-r))*add

    d[n,-2]=0; d[n,-1]=0; d[n,3*n/2+1]=0
\end{Verbatim}
}

\begin{proof}
Our problem \(see \eqref{problemd}\) is to compute
$\texttt{d[n,k]}=d^{(n)}_k+\Turing(\gamma_{n,k})$ for $0\le n\le
L-1$, $0\le k\le 3n/2$.  where $\gamma_{n,k}$ is given in
\eqref{defgammak}.  Assuming that we have computed for some even $n$
all the coefficients \texttt{d[n,k]} for $0\le k<3n/2$ with this
accuracy, we will apply then formula \eqref{case3n2k} to compute
\texttt{d[n,3n/2]}.  It is easy to see that we would not get in this
way the desired accuracy $\gamma_{n,3n/2}$. We are forced then to
compute these coefficients to a slightly better accuracy.

Define $f_x=\frac{1}{(x+1)(x+2)}$.  We substitute our problem by
\begin{equation}\label{accuracyd2}
\texttt{d[n,k]}=d^{(n)}_k+\Turing\Bigl(f_{3n/2-k}\frac{\sqrt{2\pi}}{128}\frac{(\pi
a^2/8)^{n/2}\, 2^{2k}\,\epsilon_4} {\Gamma((3n-2k+1)/2)}\Bigr).
\end{equation}
Since $f_x\le1$ for $x\ge0$ this will be more than needed.

First we show that our precision \texttt{wpd[n]} is decreasing with
$n$. This is equivalent to show that $g_n:=(8/\pi a^2
B_1^2)^{n/2}\Gamma(n-\tfrac12 )^{1/2}$ is decreasing with $n$. But by
\eqref{thirdcondition} and the facts that  $n<L$, $B_1 \ge1$ and
$L\ge2$
\begin{displaymath}
\frac{g_n}{g_{n+1}}=  a B_1\frac{\sqrt{\pi/8}}{\sqrt{n-1/2}}>
\Bigl(\frac{75
L}{2}\Bigr)^{1/2}\frac{\sqrt{\pi/8}}{\sqrt{L-3/2}}>\sqrt{18\pi}>1.
\end{displaymath}

We proceed by induction on $n$. We assume that we have computed
\texttt{d[n-1,k]} for $-2\le k\le \lfloor3(n-1)/2\rfloor+1$
satisfying \eqref{accuracyd2} and we want to compute \texttt{d[n,k]}
for $-2\le k\le \lfloor3n/2\rfloor+1$. This is true in the first run
\(for $n=1$\) since the values of $d^{(0)}_k$ for $-2\le k\le 1$ are
given exactly on line 4.

So we consider the loop starting in line 5.  Since
$\texttt{wpd[n]}\le \texttt{wpd[n-1]}$ we will have a truncation
error in each value computed previously.

The values \texttt{d[n, k]} for $k=-2$, $k=-1$ and
$k=\lfloor3n/2\rfloor+1$ are computed exactly  on line 20. So we
only need to consider the cases $0\le k\le 3n/2$.

First consider the case where $0\le k<3n/2$. In this case
$m:=3n-2k>0$. In line 10--12 we get the values
\begin{displaymath}
\texttt{m1}=\frac{1}{m}(1+\eta_1),\quad
\texttt{c1}=\frac{1}{4m}(1+\eta_1),\quad
\texttt{c2}=\frac{1-2\sigma}{2m}(1+\eta_3),\quad \texttt{c3}=-(m+1).
\end{displaymath}
By the induction hypothesis we have
$\texttt{d[n-1,k]}=d^{(n-1)}_k+\alpha_k$ where $|\alpha_k|\le
f_{3n/2-k}\gamma_{n-1,k}$ \(where $\gamma_{n,k}=0$ for $k=-2$,
$k=-1$ and $k=\lfloor 3n/2\rfloor+1$\). Then in line 14 we compute
\begin{multline*}
\texttt{d[n,k]}=-(m+1)(d^{(n-1)}_{k-2}+\alpha_{k-2})(1+\eta_5)+
\frac{1}{4m}(d^{(n-1)}_{k}+\alpha_{k})(1+\eta_5)+\\
+ \frac{1-2\sigma}{2m}(d^{(n-1)}_{k-1}+\alpha_{k-1})(1+\eta_6)
\end{multline*}
so that by \eqref{recurrenced} we have $\texttt{d[n,k]} = d^{(n)}_k
+R_k$, where
\begin{multline*}
R_k=-(m+1)(d^{(n-1)}_{k-2}+\alpha_{k-2})\eta_5-(m+1)\alpha_{k-2}+
\frac{1}{4m}(d^{(n-1)}_{k}+\alpha_{k})\eta_5+\\
+\frac{1}{4m}\alpha_k+\frac{1-2\sigma}{2m}(d^{(n-1)}_{k-1}+\alpha_{k-1})\eta_6+
\frac{1-2\sigma}{2m}\alpha_{k-1}
\end{multline*}

We put $R_k=R'_k+R''_k$ where $R'_k$ consists of the terms which
contains a factor $d^{(n)}_k$ and $R''_k$ consists of the other
terms. First we consider $R'_k$.  We will apply in this case the
bound $|d^{(n)}_k|\le D^{(n)}_k$ given in \cite{A86}*{equation (49)}. \(Observe that $D^{(n)}_k>0$ when $0\le k\le 3n/2$ and
can be taken equal to $0$ when this condition is not satisfied.\)
Then
\begin{equation}
|R'_k|\le
(m+1)D^{(n-1)}_{k-2}|\eta_5|+\frac{1}{4m}D^{(n-1)}_{k}|\eta_5|+\frac{|1-2\sigma|}{2m}
D^{(n-1)}_{k-1}|\eta_6|
\end{equation}
Substituting the value of $D^{(n)}_k$ and taking common factors we
get

\[
|R'_k|\le
A\frac{2^k}{B_1^{n-1}}\Bigl(\frac{\Gamma(n-\frac12)}{m!}\Bigr)^{1/2}
\Bigl\{\frac{\sqrt{m+1}}{4}|\eta_5|+\frac{\sqrt{(m-1)(m-2)}}{4\sqrt{m}}|\eta_5|+
\frac{|1-2\sigma|}{4\sqrt{m}}|\eta_6|\Bigr\}.
\]
\(This is true even when $m=1$ or $m=2$, these corresponds to the
cases in which $D^{(n-1)}_k=0$.\)

By our choice $\texttt{wpd[n]}\ge6$ and Lemma \ref{boundExp}
applies, so that
\begin{multline*}
|R'_k|\le
A\frac{2^{k-d}}{B_1^{n-1}}\Bigl(\frac{\Gamma(n-\frac12)}{m!}\Bigr)^{1/2}\cdot\\
\cdot
\Bigl\{1.325\sqrt{m+1}+1.325\frac{\sqrt{(m-1)(m-2)}}{\sqrt{m}}+
1.59\frac{|1-2\sigma|}{\sqrt{m}}\Bigr\}.
\end{multline*}

We want
\begin{equation}\label{parcialR}
|R'_k|\le \frac{f_{m/2}}{2}\cdot\frac{\sqrt{2\pi}}{128}\frac{(\pi
a^2/8)^{n/2}\, 2^{2k}\,\epsilon_4} {\Gamma((m+1)/2)}
\end{equation}
This is satisfied if we choose the working precision $d$ such that
\begin{displaymath}
2^d\ge \frac{256 A B_1 2^{-k}(8/\pi a^2
B_1^2)^{n/2}}{\sqrt{2\pi}\varepsilon_4}\Gamma(n-\tfrac12 )^{1/2} h(m)
\end{displaymath}
where
\[
h(m):= \frac{\Gamma((m+1)/2)}{f_{m/2}\sqrt{m!}}
\Bigl\{1.325\sqrt{m+1}+1.325\frac{\sqrt{(m-1)(m-2)}}{\sqrt{m}}+
1.59\frac{|1-2\sigma|}{\sqrt{m}}\Bigr\}.
\]
It is easy to see that $h(m)\le 17.5203+11.925|\sigma|$. Since
$B_1\le 2\sqrt{1-\log2}$ we get
\begin{displaymath}
\frac{256 B_1}{\sqrt{2\pi}}(17.5203+11.925|\sigma|)\le
2^{11}(1+|\sigma|).
\end{displaymath}
Hence have \eqref{parcialR} if we take
\begin{equation}\label{firstconditionind}
2^d\ge 2^{11}(1+|\sigma|)\frac{A}{ \varepsilon_4}\Bigl(\frac{8}{\pi
a^2 B_1^2}\Bigr)^{n/2} \Gamma(n-\tfrac12 )^{1/2}.
\end{equation}

Now we must bound $R''_k$
\begin{displaymath}
|R''_k|=\Bigl|-(m+1)\alpha_{k-2}(1+\eta_5)+\frac{1}{4m}\alpha_k(1+\eta_5)
+\frac{1-2\sigma}{2m}\alpha_{k-1}(1+\eta_6)\Bigr|.
\end{displaymath}
By the induction hypothesis we have \(since $d\ge6$, $|1+\eta_6|\le
1.1$,  observe also that $m=3n-2k$\)
\begin{multline*}
\frac{|R''_k|}{1.1}\le
(m+1)f_{(m+1)/2}\frac{\sqrt{2\pi}}{128}\frac{(\pi a^2/8)^{(n-1)/2}\,
2^{2k-4}\,\epsilon_4} {\Gamma((m+2)/2)}+\\+\frac{f_{(m-3)/2}}{4m}
\frac{\sqrt{2\pi}}{128}\frac{(\pi a^2/8)^{(n-1)/2}\,
2^{2k}\,\epsilon_4}
{\Gamma((m-2)/2)}+\\
+f_{(m-1)/2}\frac{1+2|\sigma|}{2m} \frac{\sqrt{2\pi}}{128}\frac{(\pi
a^2/8)^{(n-1)/2}\, 2^{2k-2}\,\epsilon_4} {\Gamma(m/2)}
\end{multline*}
so that
\begin{multline*}
|R''_k|\le \frac{f_{m/2}}{2}\cdot\frac{\sqrt{2\pi}}{128}\frac{(\pi
a^2/8)^{n/2}\, 2^{2k}\,\epsilon_4} {\Gamma((m+1)/2)}\cdot\\
\cdot \frac{1.1}{(\pi
a^2/8)^{1/2}}\Bigl\{\frac{(m+1)\Gamma((m+1)/2)}{2^3 \Gamma((m+2)/2)}
\frac{f_{(m+1)/2}}{f_{m/2}}+ \frac{\Gamma((m+1)/2)}{2m
\Gamma((m-2)/2)} \frac{f_{(m-3)/2}}{f_{m/2}}+\\
+ \frac{(1+2|\sigma|)\Gamma((m+1)/2)}{4m \Gamma(m/2)}
\frac{f_{(m-1)/2}}{f_{m/2}}\Bigr\}
\end{multline*}
By \eqref{thirdcondition} we have
\begin{displaymath}
m+1=3n -2k+1\le 3n+1\le 3L-2\le 3L<\frac{2a^2}{25}.
\end{displaymath}
Hence
\begin{multline*}
|R''_k|\le \frac{f_{m/2}}{2}\cdot\frac{\sqrt{2\pi}}{128}\frac{(\pi
a^2/8)^{n/2}\, 2^{2k}\,\epsilon_4} {\Gamma((m+1)/2)}
\Bigr\{1.1\Bigl(\frac{16}{25\pi}\Bigr)^{1/2}\cdot\\
\cdot\Bigl(\frac{\sqrt{m+1}\Gamma((m+1)/2)}{2^3 \Gamma((m+2)/2)}
\frac{f_{(m+1)/2}}{f_{m/2}}+ \frac{\Gamma((m+1)/2)}{2m\sqrt{m+1}
\Gamma((m-2)/2)} \frac{f_{(m-3)/2}}{f_{m/2}}+\\
+ \frac{\Gamma((m+1)/2)}{4m\sqrt{m+1} \Gamma(m/2)}
\frac{f_{(m-1)/2}}{f_{m/2}}\Bigr)+\frac{1.1}{(\pi/8)^{1/2}}
\frac{2|\sigma|\Gamma((m+1)/2)}{4am \Gamma(m/2)}
\frac{f_{(m-1)/2}}{f_{m/2}}\Bigr\}
\end{multline*}

We assume that $|\sigma|\le \frac{a}{2}$ we can substitute
$\frac{|\sigma|}{a}$ by $1/2$.  The the expression between
$\{\dots\}$ is then a   function of $m$. This function has a finite
limit when $m\to+\infty$ \($= 11/25\sqrt{2\pi}$ for the curious\)
and it not difficult to see that for $m\ge0$ de maximum of the
expression between $\{\dots\}$ is attained for $m= 1$ and it is
equal to $0.618967$. So we have
\begin{equation}\label{parcialR2}
|R''_k|< \frac{f_{m/2}}{2}\cdot\frac{\sqrt{2\pi}}{128}\frac{(\pi
a^2/8)^{n/2}\, 2^{2k}\,\epsilon_4} {\Gamma((m+1)/2)}
\end{equation}
Both \eqref{parcialR} and \eqref{parcialR2} implies that
\begin{equation}\label{firstcase}
\texttt{d[n,k]}=d^{(n)}_k+\Turing(f_{3n/2-k}\gamma_{n,k}),\qquad
0\le k<3n/2:
\end{equation}

Finally in the case $n$ even and $3n=2k$, that is when $m=0$ we just
had computed \texttt{d[n,k]} satisfying \eqref{firstcase} and with
these round numbers we operate in lines 16--19 to compute the last
\texttt{d[n,k]}.  Following the rules \eqref{rulesum}   we get
\begin{displaymath}
\texttt{d[n,k]}=\sum_{r=0}^{k-1}(-1)^{k-r}(d^{(n)}_r+\beta_r)\frac{(2k-2r)!}{(k-r)!}
\frac{(1+\eta_1)}{(1+\eta_1)}(1+\eta_1) (1+\eta_{k-r})
\end{displaymath}
where $|\beta_r|\le f_{3n/2-r}\gamma_{n,r}$. By Lemma \ref{lemma2}
we have
\begin{displaymath}
\texttt{d[n,k]}=\sum_{r=0}^{k-1}(-1)^{k-r}(d^{(n)}_r+\beta_r)\frac{(2k-2r)!}{(k-r)!}
(1+\eta_{k-r+4}).
\end{displaymath}
Hence
\begin{multline}\label{diferenceexpres}
\texttt{d[n,k]}-d^{(n)}_k=\\
=\sum_{r=0}^{k-1}(-1)^{k-r}\beta_r\frac{(2k-2r)!}{(k-r)!}
(1+\eta_{k-r+4})+\sum_{r=0}^{k-1}(-1)^{k-r}d^{(n)}_r\frac{(2k-2r)!}{(k-r)!}
\eta_{k-r+4}.
\end{multline}
Call $T_1$ and $T_2$ the above two sums. Then
\begin{displaymath}
|T_1|\le\sum_{r=0}^{k-1}\frac{\sqrt{2\pi}}{128}\frac{f_{3n/2-r}(\pi
a^2/8)^{n/2}\, 2^{2r}\,\epsilon_4}
{\Gamma((3n-2r+1)/2)}\frac{(2k-2r)!}{(k-r)!} (1+|\eta_{k-r+4}|)
\end{displaymath}
Observe that in this case $3n=2k$ so that by the formula of
duplication of Legendre
\begin{multline*}
\Gamma((3n-2r+1)/2) (k-r)! =
\Gamma((3n-2r+1)/2)\Gamma(k-r+1)=\\
=(k-r)2^{1-2k+2r}\sqrt{\pi}\Gamma(2k-2r)=
2^{2r-2k}\sqrt{\pi}(2k-2r)!
\end{multline*}
Hence
\begin{displaymath}
|T_1|\le\sum_{r=0}^{k-1}\frac{\sqrt{2}}{128}f_{3n/2-r}(\pi
a^2/8)^{n/2}\, 2^{2k}\,\varepsilon_4(1+|\eta_{k-r+4}|)
\end{displaymath}
Since $k-r+4\le k+4=\frac{3n}{2}+4\le \frac{3L+5}{2}\le 2L+2<4L^2$
and we assume that $4L^2 2^{-d}<0.1$, we have $|\eta_{k-r+4}|\le
1.06\cdot (k-r+4)\cdot 2^{-d}$. Therefore the above sum is
\begin{multline}\label{boundT1}
|T_1|\le\sum_{r=0}^{k-1} \frac{\sqrt{2}}{128}f_{3n/2-r}(\pi
a^2/8)^{n/2}\, 2^{2k}\,\varepsilon_4(1+1.06(k-r+4)2^{-d})=\\=
\frac{\sqrt{2\pi}}{128}\frac{(\pi a^2/8)^{n/2}\,
2^{2k}\,\varepsilon_4}{\sqrt{\pi}}\Bigl\{\frac{1}{2}-\frac{1}{k+2}+1.06(1+2\log
k)2^{-d}\Bigr\}.
\end{multline}
{\footnotesize To justify the above observe that
\begin{displaymath}
\sum_{r=0}^{k-1}f_{3n/2-r}=\sum_{r=0}^{k-1}\frac{1}{(k-r+1)(k-r+2)}=\sum_{j=1}^{k}
\frac{1}{(j+1)(j+2)}= \frac{1}{2}-\frac{1}{k+2}.
\end{displaymath}
and
\begin{multline*}
\sum_{r=0}^{k-1}f_{3n/2-r}(k-r+4)=\sum_{j=1}^k\frac{j+4}{(j+1)(j+2)}=
\sum_{j=1}^k\frac{2}{(j+1)(j+2)}+\sum_{j=1}^k\frac{1}{(j+1)}<\\
<1+\log(k+1)<1+2\log k.
\end{multline*}
}

To the second sum we apply Schwarz inequality
\begin{multline*}
|T_2|^2=\bigl|\sum_{r=0}^{k-1}(-1)^{k-r}d^{(n)}_r\frac{(2k-2r)!}{(k-r)!}
\eta_{k-r+4}\Bigr|^2\le\\ \le
\Bigl(\sum_{r=0}^{k-1}2^{2k-2r}(2k-2r)! (d^{(n)}_r)^2\Bigr)
\Bigl(\sum_{r=0}^{k-1}\frac{(2k-2r)!|\eta_{k-r+4}|^2}{2^{2k-2r}(k-r)!
(k-r)! }\Bigr)
\end{multline*}
Inside the second factor we find the binomial coefficient
$\binom{2k-2r}{k-r}\le 2^{2k-2r}$, so that
\begin{displaymath}
|T_2|^2\le k\{1.06 (k+4)
2^{-d}\}^2\Bigl(\sum_{r=0}^{k-1}2^{2k-2r}(2k-2r)!
(d^{(n)}_r)^2\Bigr).
\end{displaymath}
Now we apply the bound $|d^{(n)}_r|\le D^{(n)}_r$
\(\cite{A86}*{(49)}\), and also observe that $3n=2k$ hence $k\ge3$
and $k+4\le 3k$. In this way we get
\begin{multline}\label{boundT2}
|T_2|^2\le k\{1.06 (k+4) 2^{-d}\}^2 \Bigl(k 2^{2k} A^2 B_1^{-2n}
\Gamma(n+\tfrac12 )\Bigr)<\\
<11\, A^2 (8/B_1^2)^n  \Gamma(n+\tfrac12 ) k^4 2^{-2d}.
\end{multline}
By equation \eqref{boundT1} and \eqref{boundT2}  we will have
$|\texttt{d[n,k]}-d^{(n)}_k|\le
f_{0}\gamma_{n,k}=\frac{1}{2}\gamma{n,k}$ if
\begin{multline*}
\frac{\sqrt{2\pi}}{128}\frac{(\pi a^2/8)^{n/2}\,
2^{2k}\,\varepsilon_4}{\sqrt{\pi}}\Bigl\{-\frac{1}{k+2}+1.06(1+2\log
k)2^{-d}\Bigr\}+\\
+\sqrt{11} A (8/B_1^2)^{n/2} \Gamma(n+\tfrac12 )^{1/2} k^2 2^{-d}<0
\end{multline*}
or equivalently
\begin{displaymath}
\frac{2^d}{k+2}>1.06(1+2\log k)+ 128
A\sqrt{11/2}\frac{k^2\Gamma(n+\frac12)^{1/2}}{(\sqrt{\pi}\,a B_1)^n
\varepsilon_4}
\end{displaymath}
Here $3n=2k$ and $n\ge2$ so that $(k+2)(1+2\log k)\le 4n^2$ and
$(k+2)k^2\le 6n^3$. Hence we can write the above condition as
\begin{displaymath}
2^d>4.24n^2 + 768
A\sqrt{11/2}\frac{n^3\Gamma(n+\frac12)^{1/2}}{(\sqrt{\pi}\,a B_1)^n
\varepsilon_4}
\end{displaymath}
We separate this in two conditions
\begin{equation}\label{secondconditionind}
2^d>8.48n^2,\qquad  2^d> 1536
A\sqrt{11/2}\frac{n^3\Gamma(n+\frac12)^{1/2}}{(\sqrt{\pi}\,a B_1)^n
\varepsilon_4}.
\end{equation}
The second condition in \eqref{secondconditionind}  is in some way
comparable to \eqref{firstconditionind}. In fact the inequality
\begin{displaymath}
1536A\sqrt{11/2}\frac{n^3\Gamma(n+\frac12)^{1/2}}{(\sqrt{\pi}\,a
B_1)^n \varepsilon_4}\le 4\frac{A\,2^{11}}{
\varepsilon_4}\Bigl(\frac{8}{\pi a^2 B_1^2}\Bigr)^{n/2}
\Gamma(n-\tfrac12 )^{1/2}
\end{displaymath}
is equivalent to
\begin{displaymath}
\frac{1536\sqrt{11/2}B_1}{4\cdot 2^{11} }\le \frac{8^{n/2}}{n^3
\sqrt{n-1/2}}
\end{displaymath}
since $B_1\le 1$ this is true for all integers $n\ge2$.

It follows that \eqref{secondconditionind} and
\eqref{firstconditionind}  follows from
\begin{equation}
2^d>8.48n^2,\qquad 2^d > (1+|\sigma|)\frac{A\,2^{13}}{
\varepsilon_4}\Bigl(\frac{8}{\pi a^2 B_1^2}\Bigr)^{n/2}
\Gamma(n-\tfrac12 )^{1/2}.
\end{equation}
The first follows from our condition $2^d > 40 L^2$ and the second
is one of our conditions on $2^d$.
\end{proof}

\subsection{The sum of zeta.}

We need to compute $\sum_{n=1}^N n^{-s} +\Turing(\varepsilon_1)$. We
could use the analysis of the sum done in Section \ref{sumsection}
but instead an alternative solution that also will show the
possibilities of the package \texttt{mpmath} of Python. We will use
the function $\texttt{mp.\ba zetasum}(s,N)$ that computes to the
working precision $d=\texttt{wp}+10$ the numbers $n^{-s}$  and
return the exact sum of these numbers rounded to the precision initial working precision.

\begin{proposition}
Assume that $t>8\pi$.  In order to get the sum $\sum_{n=1}^N n^{-s}
+\Turing(\varepsilon_1)$  we define the working precision
\texttt{wpsum} as the least natural number $d$ such that
\begin{equation}
2^d>\frac{16|s|(N+N^{1-\sigma})\log N}{\varepsilon_1}
\end{equation}
and then follow the procedure:
\end{proposition}

\bigskip
\bigskip
{\footnotesize
\begin{Verbatim}[numbers= left,numbersep=2pt,numberblanklines=false,frame=single,
label=\fbox{\Large Computing the zetasum.}]

mp.prec = wpsum

S1 = mp._zetasum(s,N)
\end{Verbatim}
}

\begin{proof}
First the program  compute numbers \texttt{a[n] = exp(-s*log(n))}, so that
\begin{displaymath}
a[n]=e^{-s\log n (1+\eta_3)}(1+\eta_1)=n^{-s}e^{-\eta_3 s\log
n}(1+\eta_1)=n^{-s}(1+\xi).
\end{displaymath}
Then
\begin{displaymath}
|\xi|\le|e^{-\eta_3s\log n}-1|+|\eta_1 e^{-\eta_3s\log n}|.
\end{displaymath}
Let $\delta:=\varepsilon_1/2(N+N^{1-\sigma})<1/12$, since
$\varepsilon_1<1/6$. Since we assume that the working precision
$d=\texttt{wpsum}+10>5$, Lemma \ref{boundExp} applies and
\begin{displaymath}
|s\eta_3|\log n\le 3.18\cdot2^{-d}|s|\log
N\le\frac{\delta}{2}<\frac{1}{12}
\end{displaymath}
by the election of \texttt{wpsum}.

Then
\begin{displaymath}
|e^{-\eta_3s\log n}-1|\le |\eta_3|\,|s|\,\log N\sum_{k=1}^\infty
\frac{(1/12)^{k-1}}{k!}<3.4\cdot 2^{-d}|s|\log N<\frac{\delta}{2}
\end{displaymath}
and, since we assume $t>8\pi$ so that $N\ge2$ we have
\begin{displaymath}
|\eta_1 e^{-\eta_3s\log n}|\le 2^{-d} e^{1/12}\le (2^{-d}|s|\log
N)\frac{e^{1/12}}{|s|\log N}<\frac{e^{1/12}}{8\pi\log
2}\frac{\delta}{8}<\frac{\delta}{2}.
\end{displaymath}
It follows that $\texttt{a[n]}=n^{-s}(1+\xi)$ with $|\xi|< \delta$.

Now we add these $a[n]$  and round to \texttt{wpsum} the resulting
number, so that
\begin{displaymath}
\texttt{S1}=(1+\eta_1)\sum_{n=1}^N \texttt{a[n] } =(1+\eta_1)\sum_{n=1}^N n^{-s}(1+\xi)
\end{displaymath}
where $|\eta_1|<2^{\texttt{wpsum}}$  and with a different $\xi$ on
each summand, so that
\begin{multline*}
\Bigl|\texttt{S1}-\sum_{n=1}^N n^{-s}\Bigr|\le
2^{-\texttt{wpsum}}\sum_{n=1}^N
n^{-\sigma}(1+\delta)+\delta\sum_{n=1}^N n^{-\sigma}\le\\
\le (N+N^{1-\sigma})\{ 2^{-\texttt{wpsum}}(1+\delta)+\delta\} \le
\frac{\varepsilon_1}{16|s|\log N}(1+\delta)+
\frac{\varepsilon_1}{2}<\varepsilon_1
\end{multline*}
since $\delta<1/12$,  $|s|>1$ and $N\ge2$.
\end{proof}

\subsection{Computing $\texttt{S3}=(-1)^{N-1} U a^{-\sigma}+\Turing(\varepsilon_8)$.}

Let $\varepsilon_8=\varepsilon/3 A_2$ where $A_2 >|S_2|$ (see \eqref{defS3}).

\begin{proposition} Let $t>8\pi$ and $|\sigma|\le a/2$.  In order to compute $\texttt{S3}=(-1)^{N-1} U a^{-\sigma}+\Turing(\varepsilon_8)$ we define the working precision \texttt{wps3} as the least natural number $d$ such that
\begin{equation}
2^d>2^6\Bigl(1+\frac{2}{\varepsilon_8 a^{\sigma}}\Bigr)\frac{t}{2}\log\frac{t}{2\pi}.
\end{equation}
and follows the procedure:
\end{proposition}

\bigskip
\bigskip
{\footnotesize
\begin{Verbatim}[numbers= left,numbersep=2pt,numberblanklines=false,frame=single,
label=\fbox{\Large Computing  \texttt{S3}.}]

mp.prec = wps3

tpi = t/(2*pi)
arg = (t/2)*log(tpi)-(t/2)-pi/8
U = exp(-j*arg)
a = trunc_a(t)
asigma = pow(a, -sigma)
S3 = ((-1)**(N-1)) * asigma * U
\end{Verbatim}
}

\begin{proof}
We will have $\texttt{tpi}=\frac{t}{2\pi}(1+\eta_4)$. The analysis of the line 3 gives us
\begin{displaymath}
\texttt{arg} = \frac{t}{2}\log\bigl\{\frac{t}{2\pi}(1+\eta_4)\bigr\}(1+\eta_6)-
\frac{t}{2}(1+\eta_3)-\frac{\pi}{8}(1+\eta_2)= \frac{t}{2}
\log\frac{t}{2\pi}-\frac{t}{2}-\frac{\pi}{8}+\xi
\end{displaymath}
where
\begin{displaymath}
\xi= \frac{t}{2}
\eta_6\log\frac{t}{2\pi}+\frac{t}{2}(1+\eta_6)\log(1+\eta_4)-\frac{t}{2}\eta_3-\frac{\pi}{8}\eta_2
\end{displaymath}
For $|x|\le1/2$ we have $|\log(1+x)|\le 2|x|$, since we have $d\ge6$ Lemma \ref{boundExp} applies and
\begin{displaymath}
|\xi|\le \Bigl\{6.36\frac{t}{2}
\log\frac{t}{2\pi}+9.33\frac{t}{2}+3.18\frac{t}{2}+0.83\Bigr\}2^{-d}<19.71\cdot
2^{-d}\frac{t}{2}\log\frac{t}{2\pi}<\frac12
\end{displaymath}
for $t>8\pi>2\pi e$.

Then, \(since  $|\xi|\le 1/2$, we have $|e^{-i\xi}-1|< 1.3|\xi|$\)
\begin{displaymath}
\texttt{U}=Ue^{-i\xi}(1+\eta_1)=U+U(e^{-i\xi}-1)+Ue^{-i\xi}\eta_1=U+\delta
\end{displaymath}
where
\begin{displaymath}
|\delta|\le 1.3|\xi|+|\eta_1|<25.71\cdot
2^{-d}\frac{t}{2}\log\frac{t}{2\pi}
\end{displaymath}

In lines 5--6  we get
\begin{displaymath}
\texttt{a}=a(1+\eta_1), \quad\texttt{asigma}=\{a(1+\eta_1)\}^{-\sigma(1+\eta_1)}(1+\eta_1)=
a^{-\sigma}* (1+\kappa).
\end{displaymath}
where
\begin{displaymath}
\log (1+\kappa) =-\sigma\log(1+\eta_1)-\sigma\eta_1\log
a-\sigma\eta_1\log(1+\eta_1)+\log(1+\eta_1).
\end{displaymath}
Since $|\eta_1|\le 1/2$ we get
\begin{displaymath}
|\log(1+\kappa)|\le\{2|\sigma| +|\sigma|\log
a+2|\sigma|2^{-d}+2\}2^{-d}
\end{displaymath}
In Proposition \ref{propd} we have assumed that $|\sigma|\le
\frac{a}{2}=\frac{1}{\sqrt{2\pi t}}\frac{t}{2}$ also we assume that
$t>8\pi>2\pi e$, so that $a\le \frac{1}{4\pi}\frac{t}{2}\le
\frac{1}{4\pi}\frac{t}{2}\log\frac{t}{2\pi}$ It follows that
\begin{displaymath}
|\log(1+\kappa)|\le \{\frac{2}{\sqrt{2\pi t}} +\frac{1}{2\sqrt{2\pi
t}}+\frac{1}{\sqrt{2\pi t}}+\frac{4}{t}\}
2^{-d}\;\frac{t}{2}\log\frac{t}{2\pi}<0.44\cdot2^{-d}\;\frac{t}{2}\log\frac{t}{2\pi}.
\end{displaymath}
This is less than  $1/2$ so that
\begin{displaymath}
|\kappa|=|e^{\log(1+\kappa)}-1|<1.3\cdot|\log(1+\kappa)|<
2^{-d}\;\frac{t}{2}\log\frac{t}{2\pi}.
\end{displaymath}

Hence
\begin{displaymath}
\texttt{S3}=(-1)^{N-1}(U+\delta)a^{-\sigma}(1+\kappa)(1+\eta_1)=(-1)^{N-1}U a^{-\sigma}+R
\end{displaymath}
where
\begin{displaymath}
(-1)^{N-1}R= a^{-\sigma}U(\kappa+\eta_1+\eta_1\kappa)+\delta a^{-\sigma}(1+\kappa)(1+\eta_1).
\end{displaymath}
Let $T=\frac{t}{2}\log\frac{t}{2\pi}$, then
\begin{displaymath}
|R|\le \{a^{-\sigma}T+a^{-\sigma}+a^{-\sigma} T 2^{-d}+4\cdot 25.71\cdot a^{-\sigma} T \}2^{-d}\le  106 a^{-\sigma} T 2^{-d}<\varepsilon_8
\end{displaymath}
as we wanted.
\end{proof}

\subsection{Computation of $Z(t)$.}
We apply the formula $Z(t)=\Re\{e^{i\vartheta(t)}\Rzeta(\tfrac12 +it)\}$ valid for $t$ real.  We also assume that $t>16 \pi$.

\begin{proposition}
To compute $Z(t)$ for $t>16\pi$ real we determine two precisions
\texttt{wptheta} as the least natural number $d$ such that
\begin{equation}
d= \texttt{wpinitial}+1 + \texttt{mag}\Bigl\{3 \Bigl(\frac{t}{2\pi}\Bigr)^{3/2}\log\frac{t}{2\pi}\Bigr\}
\end{equation}
and
\texttt{wpz} as the least natural number $d$ such that
\begin{equation}
d= \texttt{wpinitial}+1 + \texttt{mag}\Bigl\{12 \frac{t}{2\pi}\log\frac{t}{2\pi}\Bigr\}
\end{equation}
and follow the procedure:
\end{proposition}

\bigskip
\bigskip
{\footnotesize
\begin{Verbatim}[numbers= left,numbersep=2pt,numberblanklines=false,frame=single,
label=\fbox{\Large Computing  $Z(t)$.}]

def RiemannSiegelZ(t):
    wpinitial = mp.prec
    mp.prec = 15
    --- compute wpz and wptheta ---
    mp.prec = wptheta
    theta = siegeltheta(t)
    mp.prec = wpz
    s=mpc(real=mpf('0.5'), imag = t)
    rzeta = Rzeta(s)
    z = exp(j*theta) * rzeta
    mp.prec = wpinitial
    return(2 * z.real)
\end{Verbatim}
}

\begin{proof}
We assume that in line 6 we compute $\texttt{theta} =
\vartheta(1+\eta'_1)$. We put $\eta'_1$ to indicate that this
is relative to the precision \texttt{wptheta}.
Observe that we will have
\begin{displaymath}
|e^{i\vartheta\eta'_1}-1|\le 1.3 |\vartheta|\eta'
\end{displaymath}

In line 8 we get
$\texttt{rzeta}=\Rzeta+\alpha=\Rzeta(\tfrac12 +it)+\alpha$ with
$|\alpha|\le 2^{-d}$ where $d$ is the value of \texttt{wpz}.

So the value we return is the real part of
\begin{displaymath}
e^{i\vartheta(1+\eta'_1)}(\Rzeta+\alpha)(1+\eta_2)=e^{i\vartheta}\Rzeta+R
\end{displaymath}
where
\begin{displaymath}
R=\alpha e^{i\vartheta(1+\eta'_1)}(1+\eta_2)+\eta_2\Rzeta
e^{i\vartheta(1+\eta'_1)}+e^{i\vartheta}(e^{i\vartheta
\eta'_1}-1)\Rzeta
\end{displaymath}
It can be shown that for $\sigma>0$ and $t>16\pi$ we have
$|\Rzeta(\sigma+it)|<2\sqrt{t/2\pi}$, also
$|\vartheta(t)|<T:=\frac{t}{2}\log\frac{t}{2\pi}$
Hence,
\[
|R|\le3\cdot 2^{-d}+ 2.12\,  T 2^{-d}+1.3\, T 2^{-d'}
2\sqrt{\frac{t}{2\pi}}<6\, T 2^{-d}+\frac{1}{2}
2^{-\texttt{wpinitial} }<2^{-\texttt{wpinitial}}.\qedhere
\]
\end{proof}

\subsection{Computing $\zeta(s)$. }
We shall apply the formula
\begin{equation}\label{zetaRR}
\zeta(\sigma+it)=\Rzeta(\sigma+it)+e^{-2i\vartheta(t-i(\sigma-\frac12))}
\overline{\Rzeta(1-\sigma+it)}.
\end{equation}

With the usual notations $\chi(\sigma+it)=e^{-2i\vartheta(t-i(\sigma-\frac12))}$. We need some simple bounds of $|\Rzeta(s)|$ and $|\chi(s)|$. We shall use.
\begin{proposition} For $\sigma>0$, $t>1/2$ and $|s|\ge2\pi e$
\begin{equation}\label{chimas}
|\chi(\sigma+it)|\le (\sigma^2+t^2)^{\frac{1}{4}}.
\end{equation}
\end{proposition}
\begin{proposition} For $\sigma<0$ and $t>\tfrac12 $ we have
\begin{equation}\label{chimenos}
|\chi(\sigma+it)|\le
\frac{1}{(2\pi)^{1-\sigma}}\{(1-\sigma)^2+t^2\}^{-\frac{\sigma}{2}+\frac{1}{4}}.
\end{equation}
\end{proposition}

\begin{proposition} We have
\begin{equation}\label{Rright}
|\Rzeta(\sigma+it)|\le 2\sqrt{\frac{t}{2\pi}}\qquad \sigma>0, \quad t>16\pi.
\end{equation}
\end{proposition}
\begin{proposition}
 For $\sigma<0$ and $t>16\pi$ we have
\begin{equation}\label{Rzetaleft}
|\Rzeta(\sigma+it)|\le
\frac{4t}{(2\pi)^{1-\sigma}}\{(1-\sigma)^2+t^2\}^{-\frac{\sigma}{2}+\frac{1}{4}}\qquad
\sigma<0,\quad t>16\pi.
\end{equation}
\end{proposition}
\begin{proposition}\label{TProp}
We have for complex $t$
\begin{equation}
|\vartheta(t)| \le 2|t|\log|t|,\qquad |t|>4,\quad |\Re t|\ge 1.
\end{equation}
\end{proposition}
\noindent whose proofs may be found in \cite{A92}.

\begin{proposition}
To compute $\zeta(s)+\Turing(\varepsilon)$ we shall determine
by the above Propositions bounds
\begin{displaymath}
|\Rzeta(\sigma+it)|\le M_1,\quad
|\Rzeta(1-\sigma+it)|\le M_2,\quad |\chi(\sigma+it)|\le X,
\end{displaymath}
and
\begin{displaymath}
|\vartheta(t+i(\tfrac12 -\sigma))|\le T
\end{displaymath}
and with them
three working precisions \texttt{wptheta}, \texttt{wpR} and \texttt{wpbasic}.
\begin{align*}
\texttt{wpbasic}&=\max\{6,3+\texttt{mag}(T), 2+\texttt{mag}(2.12 M_1+21.2 R_2X+1.3 M_2XT)-\\
&\hspace{7cm}-\texttt{mag}(\varepsilon)+1\},\\
\texttt{wptheta}&=\max\{4, 3+\texttt{mag}(2.7M_2X)-\texttt{mag}(\varepsilon)+1\},\\
\texttt{wpR}&=3+\texttt{mag}(1.1+2X)-\texttt{mag}(\varepsilon)+1.
\end{align*}
and with them follow the procedure:
\end{proposition}

\bigskip
\bigskip
{\footnotesize
\begin{Verbatim}[numbers= left,numbersep=2pt,numberblanklines=false,frame=single,
label=\fbox{\Large Computing  $\zeta(s)$.}]

def zeta(s):
    sigma = s.real
    t = s.imag
    wpinitial = mp.prec
    mp.prec = 53
    --- compute wptheta, wpR, wpbasic ---
    mp.prec = wptheta
    theta = siegeltheta(t-j*(sigma-mpf('0.5')))
    mp.prec = wpR
    rzeta = Rzeta(s)
    rzeta2 = conj( Rzeta(1-sigma+j* t) )
    mp.prec = wpbasic
    zv = rzeta+exp(-2*j*theta)*rzeta2
    return(zv)
\end{Verbatim}
}

\begin{proof}
Let $\eta_k$, $d$ refers to the precision \texttt{wpbasic}, $\eta'_k$ and $d'$ to \texttt{wptheta} and $\eta''_k$ and $d''$ to \texttt{wpR}. Then we will have
$\texttt{theta}=\vartheta(t-i(\sigma-\frac{1}{2}))+\alpha$ with $|\alpha|<2^{-d'}$, $\texttt{rzeta}=\Rzeta(s)+\beta$ and
$\texttt{rzeta2}=\overline{\Rzeta(1-\sigma+it)}+\gamma$ with $|\beta|$ and $|\gamma|<2^{-d''}$.

Then with $R_1:=\Rzeta(s)$, $R_1:=\overline{\Rzeta(1-\sigma+it)}$, $\vartheta:=\vartheta(t-i(\sigma-\frac{1}{2}))$  we will have
\(observe that both \texttt{rzeta2} and $e^{-2i(\vartheta+\alpha)}$ are complex numbers\)
\begin{displaymath}
\texttt{zv} =(R_1+\beta)(1+\eta_2)+e^{-2i(\vartheta+\alpha)(1+\eta_1)}
(R_2+\gamma)(1+\eta_5)
\end{displaymath}
so that $\texttt{zv} = \zeta(s)+E$ with
\begin{multline*}
E=R_1\eta_2+\beta(1+\eta_2)+e^{-2i\vartheta}
e^{-2i\vartheta\eta_1-2i\alpha(1+\eta_1)}\{(R_2+\gamma)\eta_5+\gamma\}
+\\+e^{-2i\vartheta}
(e^{-2i\vartheta\eta_1-2i\alpha(1+\eta_1)}-1)R_2
\end{multline*}
Since we choose $\texttt{wpbasic}\ge6$ we can apply Lemma \ref{boundExp} so that
\begin{multline*}
|E|\le 2.12 M_1 2^{-d}+1.1 \cdot2^{-d''}+X|e^{-2i\vartheta\eta_1-2i\alpha(1+\eta_1)}|(2M_2 5.3\cdot 2^{-d}+2^{-d''})+\\+
X
|e^{-2i\vartheta\eta_1-2i\alpha(1+\eta_1)}-1| M_2
\end{multline*}
Let $|\vartheta|\le T$. We choose $d$ and $d'$ so that   $2T2^{-d}+2.04\cdot 2^{-d'}\le1/2$, so that
\begin{displaymath}
X
|e^{-2i\vartheta\eta_1-2i\alpha(1+\eta_1)}-1| M_2\le 1.3 M_2 X(T2^{-d}+2.04\cdot 2^{-d'}).
\end{displaymath}
Also we have $|e^{-2i\vartheta\eta_1-2i\alpha(1+\eta_1)}|\le 2$. Hence
\begin{displaymath}
|E|\le (2.12 M_1+21.2 R_2X+1.3 M_2XT)2^{-d}
+(1.1+2X)2^{-d''} +2.7M_2X2^{-d'}
\end{displaymath}
And by our choice of precisions we get $|E|<\varepsilon$.
\end{proof}

In practice the computation of $\Rzeta(\sigma+it)$ and $\Rzeta(1-\sigma+it)$ is done simultaneously. The sum of zeta can be simplified to almost half using that $n^{-\sigma-it}\cdot \overline{n^{-1+\sigma-it}}=n^{-1}$.
Also the numbers $F^{(m)}(p)$ are exactly the same in the two cases. The time spent on the computation of $F^{(m)}(p)$ is an important part of the total time of computation of $\Rzeta(s)$.

\subsection{Computing $Z(w)$ for complex $w$.}

To compute $Z(\sigma+it)$ we may apply
\begin{equation}\label{zRR}
Z(t-i(\sigma-\frac12))=e^{i\vartheta(t-i(\sigma-\frac12))}\Rzeta(\sigma+it)
+e^{-i\vartheta(t-i(\sigma-\frac12))}
\overline{\Rzeta(1-\sigma+it)}
\end{equation}
that easily follows from \eqref{zetaRR}.  So to compute $Z(w)$ we put $t=\Re w$, $\sigma=\tfrac12 -\Im w$ and $s=\sigma+it$.

Then the  procedure is similar to that for computing  $\zeta(s)$. So we shall compute $\texttt{theta}=\vartheta(t-i(\sigma-\frac{1}{2}))+\alpha$ with $|\alpha|<2^{-d'}$, $\texttt{rzeta}=\Rzeta(s)+\beta$ and
$\texttt{rzeta2}=\overline{\Rzeta(1-\sigma+it)}+\gamma$ with $|\beta|$ and $|\gamma|<2^{-d''}$.  As before we put $R_1:=\Rzeta(s)$, $R_1:=\overline{\Rzeta(1-\sigma+it)}$, $\vartheta:=\vartheta(t-i(\sigma-\frac{1}{2}))$.
At the end we will have
\begin{displaymath}
\texttt{zv}=e^{i(\vartheta+\alpha)(1+\eta_1)}(R_1+\beta)(1+\eta_6)+
e^{-i\vartheta(1+\alpha)(1+\eta_1)}(R_2+\gamma)(1+\eta_6):=Z(t)+E.
\end{displaymath}
\(7 is 1 from the exponential, 1 of truncating $(R_1+\beta)$, 3 from the product, 1 from the last sum.\)
It follows that
\begin{multline*}
E=e^{i\vartheta} e^{i\vartheta\eta_1+i\alpha(\eta_1+1)}(R_1+\beta)\eta_6+
e^{i\vartheta} e^{i\vartheta\eta_1+i\alpha(\eta_1+1)}\beta+\\+
e^{i\vartheta}(e^{i\vartheta\eta_1+i\alpha(\eta_1+1)}-1)R_1
+e^{-i\vartheta} e^{i\vartheta\eta_1+i\alpha(\eta_1+1)}(R_2+\gamma)\eta_6+\\+
e^{-i\vartheta} e^{i\vartheta\eta_1+i\alpha(\eta_1+1)}\gamma+
e^{-i\vartheta}(e^{i\vartheta\eta_1+i\alpha(\eta_1+1)}-1)R_2
\end{multline*}
We assume that $d\ge6$ and Lemma \eqref{boundExp} applies.
Also we assume $|\beta|\le 2^{-d''}\le M_1$,
$|\gamma|\le M_2$, and we choose $d$ and $d'$ so that   $2T2^{-d}+2.04\cdot 2^{-d'}\le1/2$. That is such that $|\vartheta\eta_1+\alpha(\eta_1+1)|\le0.5$. With these assumptions we will have
\begin{multline*}
|E|\le2\, X^{1/2} 2(M_1+M_2)\cdot 6.36 \cdot 2^{-d}+2X^{1/2}2\cdot2^{-d''}+\\+X^{1/2}2(T2^{-d}+1.02\cdot 2^{-d'})(M_1+M_2)\le\\
\le X^{1/2}(M_1+M_2)(26 + 2T)2^{-d}+2.04
X^{1/2}(M_1+M_2)2^{-d'}+4 X^{1/2}2^{-d''}
\end{multline*}
and this will be less than $\varepsilon$ if we take the three precisions as
\begin{align}
\texttt{wpbasic}&=\max\{6,3+\texttt{mag}(T), \notag\\
&\hspace{2cm}\texttt{mag}(X^{1/2}(M_1+M_2)(26 + 2T))-\texttt{mag}(\varepsilon)+3\}, \\
\texttt{wptheta}&=\max\{4,\texttt{mag}(2.04\,X^{1/2}(M_1+M_2))-\texttt{mag}(\varepsilon)+3\},\\
\texttt{wpR}&=\texttt{mag}(4 X^{1/2})-\texttt{mag}(\varepsilon)+3
\end{align}

In practice we use only one program to compute the two functions  $\zeta(s)$ and $Z(s)$, changing only the precisions and the end of the computation.

\end{document}